\documentclass{article}
\usepackage[utf8]{inputenc}
\usepackage{epsfig,amssymb}
\usepackage[normalem]{ulem}
\usepackage[margin=0.95in]{geometry}
\usepackage{hyperref}
\usepackage{graphicx}
\usepackage{epstopdf}
\usepackage{amsmath}
\usepackage{bm}
\usepackage{setspace}

\usepackage{mathrsfs}
\usepackage{subcaption}
\usepackage[dvipsnames]{xcolor}
\usepackage{mathtools}
\usepackage{multirow}
\usepackage{lipsum}
\usepackage{comment}
\usepackage{enumerate}
\usepackage[affil-it]{authblk}
\usepackage{cite} 
\usepackage{leftidx}
\usepackage{relsize}
\usepackage{braket}
\usepackage{relsize}
\usepackage{float}
\DeclarePairedDelimiter{\ceil}{\lceil}{\rceil}
\usepackage[autostyle]{csquotes}

\usepackage{fontenc}
\usepackage{comment}
\usepackage{soul}
\usepackage[right]{lineno}
\newcommand{\FS}[1]{{\textcolor{red}{\bf[FS: #1]}}}

\newcommand{\SA}[1]{{\textcolor{olive}{\bf[SS: #1]}}}
\newcommand{\SAS}[1]{{\textcolor{blue}{#1}}}
\newcommand{\etal}{\textit{et al}.~}

\title{Geometrically Nonlinear Response of a Fractional-Order Nonlocal Model of Elasticity}
\author[1]{Sai Sidhardh}
\affil[1]{School of Mechanical Engineering, Ray W. Herrick Laboratories, Purdue University, West Lafayette, IN 47907}
\author[1]{Sansit Patnaik}
\author[1]{Fabio Semperlotti$^\dagger$}

\begin{document}
\date{}
\maketitle

\begin{abstract}
This study presents the analytical and finite element formulation of a geometrically nonlinear and fractional-order nonlocal model of an Euler-Bernoulli beam. The finite nonlocal strains in the Euler-Bernoulli beam are obtained from a frame-invariant and dimensionally consistent fractional-order (nonlocal) continuum formulation.
The finite fractional strain theory provides a positive definite formulation that results in a mathematically well-posed formulation which is consistent across loading and boundary conditions.
The governing equations and the corresponding boundary conditions of the geometrically nonlinear and nonlocal Euler-Bernoulli beam are obtained using variational principles. Further, a nonlinear finite element model for the fractional-order system is developed in order to achieve the numerical solution of the integro-differential nonlinear governing equations. Following a thorough validation with benchmark problems, the fractional finite element model (f-FEM) is used to study the geometrically nonlinear response of a nonlocal beam subject to various loading and boundary conditions. Although presented in the context of a 1D beam, this nonlinear f-FEM formulation can be extended to higher dimensional fractional-order boundary value problems.
\end{abstract}
$^\dagger$ All correspondence should be addressed to: \textit{fsemperl@purdue.edu}\\
\section{Introduction}
\label{sec: Introduction}

Several theoretical and experimental studies have shown that nonlinear and scale-dependent effects are prominent in the response of several structures ranging from sandwich composites, to layered and porous media, to random and fractal media, to media with damage and cracks, to biological materials like tissues and bones \cite{szabo1994time,adbone,fellah2004verification,stulov2016frequency}. In all these classes of problems, the coexistence of different material spatial scales renders the response fully nonlocal \cite{kroner1967elasticity,eringen1972nonlocal}. In many applications involving slender structures, large and rapidly varying loads typically induce nonlinear response. These nonlinear size-dependent effects are particularly prominent in microstructures and nanostructures that have far-reaching applications in medical devices, atomic devices, micro/nano-electromechanical devices and sensors \cite{sumelka2015fractional,civalek2016simple,rahimi2017linear}. These devices are primarily made from a combination of slender structures such as beams, plates and shells. Clearly, accurate modeling of the nonlocal and nonlinear effects in the response of these structures is paramount in many engineering as well as biomedical applications.
However, it appears that there is a limited amount of studies focusing on the modeling of geometrically nonlinear and nonlocal slender structures. In the following, we briefly review the main characteristics of these studies as well as discuss key limitations.

Nonlocal continuum theories use different strategies to enrich the classical governing equations, that describe the response of the medium at a point, with information on the behavior of distant points contained in a prescribed area of influence. Depending on the field in which this nonlocal description of the medium is applied, the same concept can be referred to with different terminologies such as action at a distance, long range interactions, or scale effects.
The key concept at the basis of nonlocal continuum theories is the horizon of influence or horizon of nonlocality which describes areas of the medium where distant points can still influence one another via long range cohesive forces \cite{kroner1967elasticity,eringen1972nonlocal}. 
Several theories based on integral methods, gradient methods, and very recently, the peridynamic approach have been developed to capture these long range interactions and analyze their effect on the response of either media or structures.
Gradient elasticity theories \cite{peerlings2001critical,aifantis2003update,guha2015review,sidhardh2018element,sidhardh2019size} account for the nonlocal behavior by introducing strain gradient dependent terms in the stress-strain constitutive law.
Integral methods \cite{polizzotto2001nonlocal,bavzant2002nonlocal,sidhardh2018effect} define constitutive laws in the form of convolution integrals between the strain and the spatially dependent elastic properties over the horizon of nonlocality.

As emphasized earlier, several applications involving nonlocal slender structures also experience geometrical nonlinearities.
Nonlocal effects on the nonlinear transverse response of shear deformable beams were modeled in \cite{yang2010nonlinear,srividhya2018nonlocal} using von K\'arm\'an nonlinearity and Eringen's nonlocal elasticity. Reddy \etal \cite{reddy2010nonlocal} modeled a nonlocal plate starting from the differential constitutive relations of Eringen and the von K\'arm\'an nonlinear strains. The differential model of Eringen's nonlocal elasticity was also used in \cite{emam2013general,lembo2016nonlinear} to analyse nonlocal and nonlinear beams. Similar models for beams and plates have also been developed in \cite{srinivasa2013model,reddy2014eringen}.
Several studies were also carried out to simulate the multiphysics response of beams such as, for example,  the nonlinear static and dynamic response of a capacitive electrostatic nanoactuator \cite{najar2015nonlinear} or the effect of nonlocal elasticity over the electroelastic response of flexoelectric nanobeams \cite{sidhardh2018effect}. 

Although the above nonlinear studies based on classical approaches have been able to capture key aspects of nonlocal elasticity, they are subject to important drawbacks. As shown in \cite{reddy2007nonlocal,challamel2008small}, the integral nonlocal methods do not allow defining a self-adjoint quadratic potential energy which results in paradoxical predictions such as hardening and absence of nonlocal effects for certain combinations of boundary conditions \cite{challamel2008small,khodabakhshi2015unified}. Note that hardening effects \textit{per se} are not inaccurate or paradoxical from a physical perspective and have been observed in several micro- and nanostructures \cite{romano2017constitutive}. However, hardening effects are not expected in classical integral approaches. This specific inaccuracy is the result of the ill-posed nature of the integral governing equations \cite{romano2017constitutive}. Also, given the complex nature of the nonlocal and nonlinear differential equations, research on the corresponding numerical methods is limited so far. It appears that only quadrature based methods \cite{najar2015nonlinear,ansari2015size} have been developed for the solution of nonlinear and nonlocal boundary value problems (BVPs).

In recent years, fractional calculus has emerged as a powerful mathematical tool to model a variety of nonlocal and multiscale phenomena. Fractional derivatives, which are a differ-integral class of operators, are intrinsically multiscale and provide a natural way to account for nonlocal effects. In fact, the order of the fractional derivative dictates the shape of the influence function (kernel of the fractional derivative) while its domain of integration defines the horizon of influence, that is the distance beyond which information is no longer accounted for in the derivative. It follows that time fractional operators enable memory effects (i.e. the response of a system is a function of its past history) while space fractional operators account for nonlocal and scale effects.
The unique set of properties of fractional operators have determined a surge of interest in the exploration of possible applications. Among the areas that have rapidly developed, we find the modeling of viscoelastic materials \cite{bagley1983theoretical,chatterjee2005statistical,patnaik2020application}, transport processes in complex media \cite{mainardi1996fractional,buonocore2018occurrence,hollkamp2019analysis}, and model-order reduction of lumped parameter systems \cite{hollkamp2018model,hollkamp2020application}.
Given the multiscale nature of fractional operators, fractional calculus has also found several applications in nonlocal elasticity. 
Space-fractional derivatives have been used to formulate nonlocal constitutive laws \cite{carpinteri2011fractional,sumelka2014thermoelasticity,sumelka2016fractional,hollkamp2019analysis} as well as to account for microscopic interaction forces \cite{laskin2006nonlinear,di2008long,cottone2009fractional}. 
Very recently fractional-order theories have been extended to model and analyze the static response, buckling characteristics, as well as the dynamic response of nonlocal beams \cite{sumelka2015fractional,alotta2014finite,rahimi2017linear,alotta2017finite,patnaik2019FEM,szajek2020selected}. However, all these studies focused on the linear geometric analysis of nonlocal structures. 

Previous works conducted on the development of nonlocal continuum theories based on fractional calculus have highlighted that the differ-integral nature of the fractional operators allows them to combine the strengths of both gradient and integral based methods while at the same time addressing a few important shortcomings of these integer-order formulations which have been highlighted earlier \cite{carpinteri2011fractional,sumelka2014thermoelasticity,sumelka2016fractional,patnaik2019FEM,patnaik2020fractional}. We highlight here that, while the fractional continuum formulation developed in this study builds upon \cite{sumelka2014thermoelasticity,sumelka2016fractional}, in our study we  have extended and generalized the formulation in order to account for asymmetric horizon conditions that typically occur near geometric/material boundaries or discontinuities. Further, note that while \cite{sumelka2016fractional} focused primarily on anisotropic nonlocality, our model is for isotropic materials with possible asymmetric horizon conditions. 
In this study, we address the shortcomings of the classical nonlocal formulations by developing a consistent nonlinear and fractional-order nonlocal beam model. Further, we develop a nonlinear fractional-order finite element method (f-FEM) to obtain the numerical solution of nonlinear and nonlocal BVPs.

The overall goal of this study is two-fold. First, we derive the governing equations for the geometrically nonlinear and nonlocal beam in strong form using variational principles. We build upon the fractional-order nonlocal continuum model proposed in \cite{patnaik2019FEM} to develop a fractional-order constitutive relation for a geometrically nonlinear Euler-Bernoulli beam. The model proposed in \cite{patnaik2019FEM} was shown to be self-adjoint, positive definite and well-posed, thus overcoming one of the important limitations mentioned above.
Further, note that the fractional-order nonlocal continuum model used in \cite{patnaik2019FEM} is frame-invariant, capable of dealing with asymmetric horizons, and it admits integer-order boundary conditions. The thermodynamic consistency of the same formulation has also been established in \cite{sidhardh2020thermoelastic}.
Second, we formulate a fully consistent and highly accurate nonlinear f-FEM to numerically investigate the response of the geometrically nonlinear fractional-order nonlocal beams. Although several FE formulations for fractional-order equations have been proposed in the literature, they are based on Galerkin or Petrov-Galerkin methods that are capable of solving only linear hyperbolic and parabolic differential equations \cite{zheng2010note,jin2016petrov}. We develop a Ritz FEM that is capable of obtaining the numerical solution of the nonlinear fractional-order elliptic boundary value problem that describes the static response of the fractional-order nonlocal and nonlinear beam. Further, by using the f-FEM we show that, independently from the boundary conditions, the fractional-order theory predicts a consistent softening behavior for the geometrically nonlinear fractional-order beam as the degree of nonlocality increases. This is a direct result of the well-posed nature of the fractional-order formulation which accepts a unique solution \cite{patnaik2019FEM}. With these results, we explain the paradoxical predictions of hardening and the absence of nonlocal effects for certain combinations of boundary conditions, as predicted by classical integral approaches to nonlocal elasticity \cite{challamel2008small,khodabakhshi2015unified}.

The remainder of the paper is structured as follows: we begin with the development of the nonlinear fractional-order beam theory based on the fractional-order von K\'arm\'an strain-displacement relations. Next, we derive the governing equations of the nonlinear fractional-order beam in strong form using variational principles. Then, we derive a strategy for obtaining the numerical solution to the nonlinear beam governing equation using fractional-order FEM. Finally, we validate the nonlinear f-FEM, establish its convergence, and then use it to analyze the effect of the fractional-order nonlocality on the geometrically nonlinear static response of the beam under different loading conditions.

\section{Nonlinear and nonlocal fractional constitutive model}
\label{sec: eb_model}
In this section, the geometrically nonlinear response of a nonlocal Euler-Bernoulli beam is modeled using the fractional-order nonlocal continuum formulation from \cite{patnaik2019FEM}. Starting from this continuum theory, the nonlinear constitutive relations for the nonlocal Euler-Bernoulli beam are derived. Finally, the governing differential equations and the corresponding boundary conditions for the fractional-order beam theory are derived using variational principles.
\subsection{Fundamentals of the fractional-order nonlocal continuum formulation}
\label{ssec: FCM}
\begin{figure}[h]
	\centering
	\includegraphics[scale=0.6]{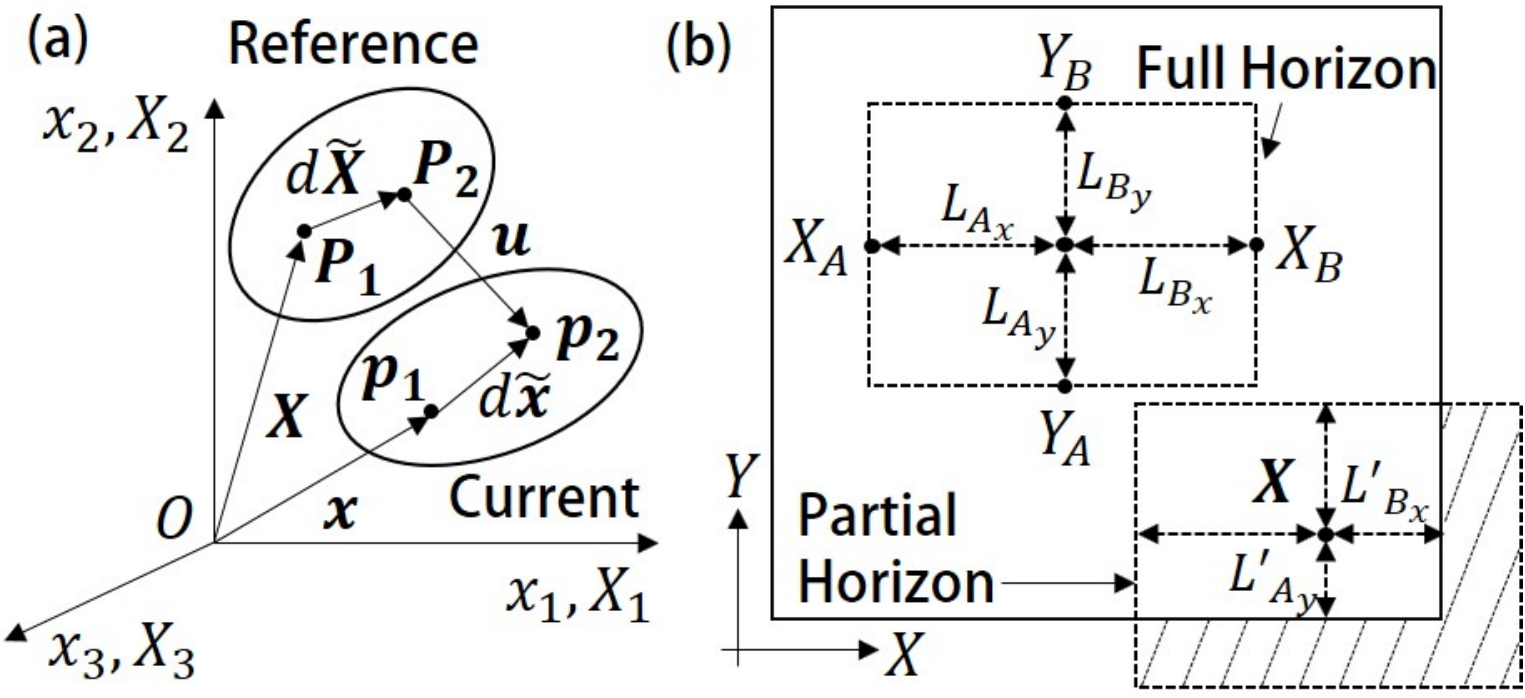}
	\caption{\label{fig1} (a) Schematic indicating the infinitesimal material $\mathrm{d}\tilde{\textbf{X}}$ and spatial $\mathrm{d}\tilde{\textbf{x}}$ line elements in the nonlocal medium under the displacement field $\textbf{u}$. (b) Horizon of nonlocality and length scales at different material points. The current nonlocal model can account for an asymmetric horizon condition that occurs at points ${\textbf{X}}$ close to a boundary or interface. The schematic is adapted from \cite{patnaik2019FEM}.}
\end{figure}

In analogy with the traditional approach to continuum modeling, we perform the deformation analysis of a nonlocal solid by introducing two stationary configurations, namely, the reference (undeformed) and the current (deformed) configurations. The motion of the body from the reference configuration (denoted as $\textbf{X}$) to the current configuration (denoted as $\textbf{x}$) is assumed as:
\begin{equation}
\label{eq: motion_eq}
\textbf{x}=\bm{\Psi}(\textbf{X})
\end{equation}
such that $\bm{\Psi}(\textbf{X})$ is a bijective mapping operation.
The relative position vector of two point particles located at $P_1$ and $P_2$ in the reference configuration of the nonlocal medium is denoted by ${\mathrm{d}\tilde{\textbf{X}}}$ (see Fig.~(\ref{fig1})). After deformation due to the motion $\bm{\Psi}(\textbf{X})$, the particles move to new positions $p_1$ and $p_2$, such that the relative position vector between them is ${\mathrm{d}\tilde{\textbf{x}}}$. It appears that $\mathrm{d}\tilde{\textbf{X}}$ and $\mathrm{d}\tilde{\textbf{x}}$ are the material and spatial differential line elements in the nonlocal medium, conceptually analogous to the classical differential line elements $\mathrm{d}{\textbf{X}}$ and $\mathrm{d}{\textbf{x}}$.

The mapping operation described in Eq.~(\ref{eq: motion_eq}) is used to model the differential line elements in the nonlocal solid using fractional calculus in order to introduce nonlocality into the system \cite{patnaik2019FEM}. More specifically, the differential line elements of the nonlocal medium are modeled by imposing a fractional-order transformation on the classical differential line elements as follows:
\begin{subequations}
\label{eq: Fractional_F}
\begin{equation}
\label{eq: Fractional_Fa}
\mathrm{d}\tilde{\textbf{x}}=\big[D^{\alpha}_\textbf{X}\bm{\Psi}(\textbf{X})\big]\mathrm{d} \textbf{X}=\big[\tilde{\textbf{F}}_{X}(\textbf{X})\big]\mathrm{d}\textbf{X}
\end{equation}
\begin{equation}
\mathrm{d}\tilde{\textbf{X}}=\big[D^{\alpha}_\textbf{x}\bm{\Psi}^{-1}(\textbf{x})\big]\mathrm{d} \textbf{x}=\big[\tilde{\textbf{F}}_{x}(\textbf{x})\big]\mathrm{d}\textbf{x}
\end{equation}	
\end{subequations}
where $D^{\alpha}_\square\square$ is a space-fractional derivative whose details will be presented below.
Given the differ-integral nature of the space-fractional derivative, the differential line elements $\mathrm{d}\tilde{\textbf{X}}$ and $\mathrm{d}\tilde{\textbf{x}}$ have a nonlocal character.
Using the definitions for $\mathrm{d}\tilde{\textbf{X}}$ and $\mathrm{d}\tilde{\textbf{x}}$, the fractional deformation gradient tensor $\overset{\alpha}{\textbf{F}}$ with respect to the nonlocal coordinates is obtained in \cite{patnaik2019FEM} as:
\begin{equation}
\label{eq: Fractional_F_net}
\frac{\mathrm{d}\tilde{\textbf{x}}}{\mathrm{d}\tilde{\textbf{X}}}=\mathop{\textbf{F}}^{\alpha}=\tilde{\textbf{F}}_{X}\textbf{F}^{-1}\tilde{\textbf{F}}_{x}^{-1}
\end{equation}
where $\textbf{F}$ is the classical deformation gradient tensor given as $\textbf{F}=\mathrm{d}\textbf{x}/\mathrm{d}\textbf{X}$ given in local and integer order form.

The space-fractional derivative $D^{\alpha}_\textbf{X}\bm{\Psi}(\textbf{X},t)$ is taken here according to a Riesz-Caputo (RC) definition with order $\alpha\in(0,1)$ defined on the interval $\textbf{X} \in (\textbf{X}_A,\textbf{X}_B) \subseteq \mathbb{R}^3 $, and given by:
\begin{equation}
\label{eq: frac_der_def}
	D^{\alpha}_\textbf{X}\bm{\Psi}(\textbf{X},t)=\frac{1}{2}\Gamma(2-\alpha)\big[\textbf{L}_{A}^{\alpha-1}~ {}^C_{\textbf{X}_{A}}D^{\alpha}_{\textbf{X}} \bm{\Psi}(\textbf{X},t) - \textbf{L}_{B}^{\alpha-1}~ {}^C_{\textbf{X}}D^{\alpha}_{\textbf{X}_{B}}\bm{\Psi}(\textbf{X},t)\big]
\end{equation}
where $\Gamma(\cdot)$ is the Gamma function, and ${}^C_{\textbf{X}_{A}}D^{\alpha}_{\textbf{X}}\bm{\Psi}$ and ${}^C_{\textbf{X}}D^{\alpha}_{\textbf{X}_{B}}\bm{\Psi}$ are the left- and right-handed Caputo derivatives of $\bm{\Psi}$, respectively. The length scales introduced in the above definition are given as $\textbf{L}_A=\textbf{X}-\textbf{X}_A$ and $\textbf{L}_B=\textbf{X}_B-\textbf{X}$ (see Fig.~(\ref{fig1})). Before proceeding, we highlight certain implications of the above definition of the fractional-order derivative.
The interval of the fractional derivative $(\textbf{X}_A,\textbf{X}_B)$ defines the horizon of nonlocality (also called attenuation range in classical nonlocal elasticity) which is schematically shown in Fig.~(\ref{fig1}) for a generic point $\textbf{X}\in\mathbb{R}^2$. The length scale parameters $L_{A_j}^{\alpha-1}$ and $L_{B_j}^{\alpha-1}$ ensure the dimensional consistency of the deformation gradient tensor, and along with the term $\frac{1}{2}\Gamma(2-\alpha)$ ensure the frame invariance of the constitutive relations. As shown in Appendix 1, the deformation gradient tensor introduced via Eq.~(\ref{eq: Fractional_F_net}) enables an efficient and accurate treatment of the frame invariance and completeness of kernel in the presence of asymmetric horizons, material boundaries, and interfaces (Fig.~(\ref{fig1})). More specifically, the use of the different length scales $L_{A_j}^{\alpha-1}$ and $L_{B_j}^{\alpha-1}$ in Eq.~(\ref{eq: frac_der_def}) allows dealing with possible asymmetries in the horizon of nonlocality (e.g. resulting from a truncation of the horizon, when approaching a boundary point or an interface). We emphasize that this should not be confused with the general concept of anisotropy wherein the length-scales would potentially be different at each point in the nonlocal continuum (see, for example, \cite{sumelka2016fractional}).

The differential line elements are used to define strain measures in the nonlocal medium. Here, we provide only the definition of strain measures in Lagrangian coordinates. In analogy with the classical strain measures, the fractional strain can be defined using the differential line elements as $d\tilde{\textbf{x}}d\tilde{\textbf{x}}-d\tilde{\textbf{X}}d\tilde{\textbf{X}}$. Using Eq.~(\ref{eq: Fractional_F_net}), the Lagrangian strain tensor in the nonlocal medium is obtained as:
\begin{equation}
\label{eq: strain_tensors}
	\mathop{\textbf{E}}^{\alpha}=\frac{1}{2}(\mathop{\textbf{F}}^{\alpha}{}^T\mathop{\textbf{F}}^{\alpha}-\textbf{I})
\end{equation}
where $\textbf{I}$ is the second-order identity tensor. 
Using kinematic position-displacement relations, the expressions of the strains can be obtained in terms of the displacement gradients. The fractional displacement gradient is obtained using the definition of the fractional deformation gradient tensor from Eq.~(\ref{eq: motion_eq}) and the displacement field $\textbf{U}(\textbf{X})=\textbf{x}(\textbf{X})-\textbf{X}$ as:
\begin{equation}
\label{eq: disp_grad_tesnor_Lag}
\nabla^\alpha {\textbf{U}}_X=\tilde{\textbf{F}}_{X}-\textbf{I}
\end{equation} 
The fractional gradient denoted by $\nabla^\alpha\textbf{U}_X$ is given as $\nabla^\alpha\textbf{U}_{X_{ij}}=D^{\alpha}_{X_j}U_i$. 
Using the nonlocal strain defined in Eq.~(\ref{eq: strain_tensors}) and the fractional deformation gradient tensor $\overset{\alpha}{\textbf{F}}$ given in Eq.~(\ref{eq: Fractional_F_net}) together with Eq.~(\ref{eq: disp_grad_tesnor_Lag}), the relationship between the strain and displacement gradient tensors is obtained as:
\begin{equation}
\label{eq: finite_fractional_strain}
\mathop{\textbf{E}}^{\alpha}=\frac{1}{2}\bigr(\nabla^\alpha {\textbf{U}}_X+\nabla^\alpha {\textbf{U}}_X^{T}+\nabla^\alpha {\textbf{U}}_X^{T}\nabla^\alpha {\textbf{U}}_X\bigl)
\end{equation}
From the strain-displacement relations, assuming moderate rotations but small strains, the fractional-order von-K\'arm\'an strain-displacement relations are defined analogously to classical elasticity theories \cite{reddy2014introduction}. The expressions for the full 3D nonlinear von-K\'arm\'an strain-displacement relations are:
\begin{equation}
\label{eq: vonkarman_strain}
    \tilde{\epsilon}_{ij}={\frac{1}{2}\left(D^{\alpha}_{X_j}u_i+D^{\alpha}_{X_i}u_j\right)}+{\frac{1}{2}\left(D^{\alpha}_{X_i}u_3~D^{\alpha}_{X_j}u_3\right)},~~~i,j=1,2
\end{equation}
where $u_3(\textbf{x})$ is the transverse displacement field and $u_k(\textbf{x})$ $(k=1,2)$ are the in-plane variables. The transverse strains (normal and shear) are simply the linearized forms of the respective expressions derived from Eq.~(\ref{eq: finite_fractional_strain}). Complete expressions for the components of the fractional-order von-K\'arm\'an strains can be written from the above equation to be:
\begin{subequations}
\begin{equation}
    \tilde{\epsilon}_{11}=D_{x_1}^\alpha u_1+\frac{1}{2}(D_{x_1}^\alpha u_3)^2
\end{equation}
\begin{equation}
    \tilde{\epsilon}_{22}=D_{x_2}^\alpha u_2+\frac{1}{2}(D_{x_2}^\alpha u_3)^2
\end{equation}
\begin{equation}
    \tilde{\epsilon}_{12}=\tilde{\epsilon}_{21}=\frac{1}{2}\left(D_{x_1}^\alpha u_2+D_{x_2}^\alpha u_1\right)+\frac{1}{2}(D_{x_1}^\alpha u_3~D_{x_2}^\alpha u_3)
\end{equation}
\begin{equation}
    \tilde{\epsilon}_{33}=D_{x_3}^\alpha u_3
\end{equation}
\begin{equation}
\label{eq: trans_e13}
    \tilde{\epsilon}_{13}=\tilde{\epsilon}_{31}=\frac{1}{2}\left(D_{x_1}^\alpha u_3+D_{x_3}^\alpha u_1\right)
\end{equation}
\begin{equation}
    \tilde{\epsilon}_{23}=\tilde{\epsilon}_{32}=\frac{1}{2}\left(D_{x_2}^\alpha u_3+D_{x_3}^\alpha u_2\right)
\end{equation}
\end{subequations}
which agree with the integer-order counterparts given in \cite{reddy2014introduction} for $\alpha=1.0$.

In this study, we have analyzed the effect of the above nonlinear strain-displacement relations on the response of an isotropic solid. The stress in this isotropic solid is given as:
\begin{equation}\label{eq: constt_isot}
    \tilde{\sigma}_{ij}=C_{ijkl}\tilde{\epsilon}_{kl}
\end{equation}
where $C_{ijkl}$ is the constitutive matrix for the isotropic solid and $\tilde{\epsilon}_{kl}$ is the fractional-order strain. We emphasize that the stress defined through the above equation is nonlocal in nature. This follows from the fractional-order definition of the deformation gradient tensor which is then reflected in the nonlocal strain as evident from Eq.~(\ref{eq: strain_tensors}). The thermodynamic consistency of the above constitutive relations has been addressed in detail in \cite{sidhardh2020thermoelastic}. As expected, classical continuum mechanics relations are recovered when the order of the fractional derivative is set as $\alpha=1$.


\subsection{Nonlinear constitutive relations}
\begin{figure}[ht!]
    \centering
    \includegraphics[width=0.6\textwidth]{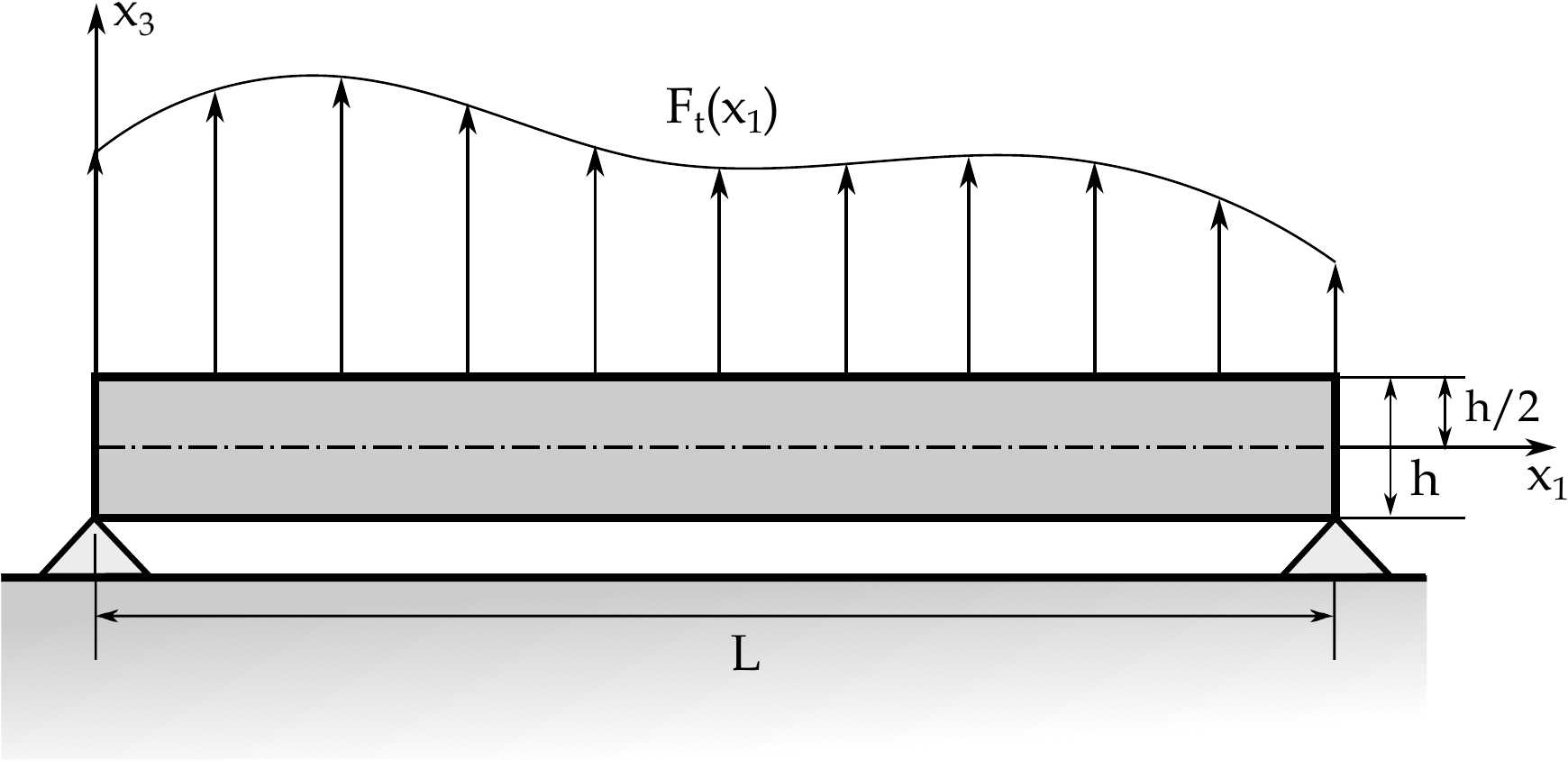}
    \caption{Schematic of an elastic beam subject to a distributed transverse load $F_t(x_1)$.}
    \label{fig: schematic_beam}
\end{figure}

We use the fractional-order nonlinear strain-displacement formulation presented above to derive the constitutive model for geometrically nonlinear nonlocal response of the beam. A schematic of a beam subject to transverse loads is shown in Fig.~(\ref{fig: schematic_beam}).
The Cartesian coordinates are chosen such that $x_3=\pm h/2$ coincide with the top and bottom surfaces of the beam, and $x_1=0$ and $x_1=L$ correspond to the two ends of the beam in the longitudinal direction. The origin of the reference frame is chosen at the left-end of the beam and $x_3=0$ is coincident with the mid-plane.

\begin{figure}[h!]
    \centering
    \includegraphics[width=0.8\textwidth]{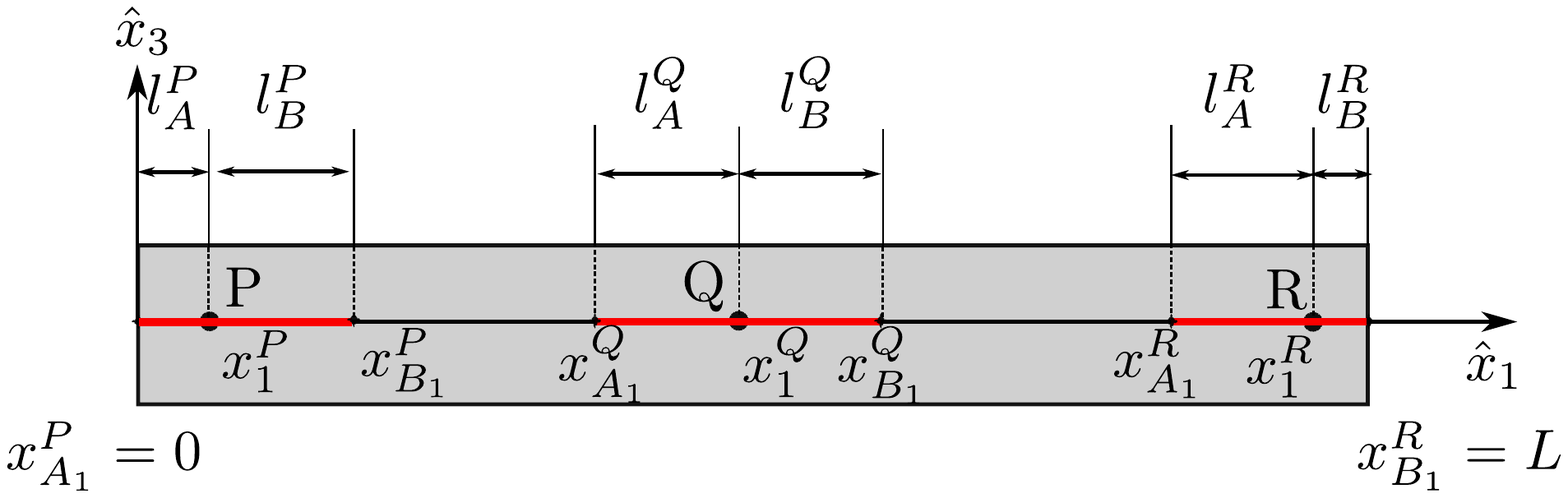}
    \caption{A schematic presenting the variable nature of the length scale corresponding to the horizon of nonlocality for different points along the length of the beam. The length scale on the left side of the point P, $l_A^P$ is truncated such that $l^P_A<l^Q_A$. Similarly, $l^R_B<l^Q_B$.}
    \label{fig: beam_horizon}
\end{figure}

In the above defined coordinate system, the axial and transverse components of the displacement field $\textbf{u}(x_1,x_3)$, denoted by $u_1(x_1,x_3)$ and $u_3(x_1,x_3)$ at any spatial location $\textbf{x}(x_1,x_3)$, are related to the mid-plane displacements according to the Euler-Bernoulli assumptions:
\begin{subequations}
    \label{eq: euler_bern_beam}
    \begin{equation}
        u_1(x_1,x_3)=u_0(x_1)-x_3\left[\frac{\mathrm{d}w_0(x_1)}{\mathrm{d}x_1}\right]
    \end{equation}
    \begin{equation}
        u_3(x_1,x_3)=w_0(x_1)
    \end{equation}
\end{subequations}
where $u_0(x_1)$ and $w_0(x_1)$ are the axial and transverse displacements of a point on the mid-plane, respectively. The fractional-order nonlinear axial strain is obtained from the fractional-order von-K\'arm\'an strain-displacement relations defined in Eq. \eqref{eq: vonkarman_strain}. For the above displacement field, the non-zero strain $\tilde{\epsilon}_{11}(x_1,x_3)$ is obtained as:
\begin{equation}
    \label{eq: vonkarman_strain_beam}
    \tilde{\epsilon}_{11}(x_1,x_3)=D_{x_1}^\alpha u_1(x_1,x_3)+\frac{1}{2}(D_{x_1}^\alpha u_3(x_1))^2
\end{equation}
where $u_1$ and $u_3$ denote the displacement fields in the $\textbf{x}_1$ and $\textbf{x}_3$ directions. For the Euler-Bernoulli displacement field assumptions in Eq. \eqref{eq: euler_bern_beam}, we obtain:
\begin{equation}
\label{eq: von_karman_ebt}
    \tilde{\epsilon}_{11}(x_1,x_3)=\tilde{\epsilon}_0(x_1)+x_3\tilde{\kappa}(x_1)
\end{equation}
where the fractional-order axial strain $\tilde{\epsilon}_0(x_1)$ and bending strain $\tilde{\kappa}(x_1)$ are given as:
\begin{subequations}
\begin{equation}
    \tilde{\epsilon}_0(x_1)=D_{x_1}^{\alpha} u_0(x_1)+\frac{1}{2}\left[D_{x_1}^{\alpha} w_0(x_1)\right]^2
\end{equation}
\begin{equation}
    \tilde{\kappa}(x_1)=-D_{x_1}^{\alpha}\left[\frac{\mathrm{d}w_0(x_1)}{\mathrm{d}x_1}\right]
\end{equation}
\end{subequations}
In the above equation, $D^\alpha_{x_1}(\cdot)$ denotes the fractional-order RC derivative given in Eq.~(\ref{eq: frac_der_def}) such that $\alpha\in(0,1)$. $\textbf{x}_A (x_{A_1},0)$ and $\textbf{x}_B(x_{B_1},0)$ are the terminals of the left- and right-handed Caputo derivatives within the RC derivative. Consequently, the domain $(\textbf{x}_A,\textbf{x}_B)$ is the horizon of nonlocal interactions for $\textbf{x} (x_1,0)$, along the mid-plane. Further, $l_A=x_1-x_{A_1}$ and $l_B=x_{B_1}-x_1$ are the length scales to the left and right hand sides of the point $\textbf{x}(x_1,0)$ along the direction $\hat{\textbf{x}}_1$. As emphasized in \S\ref{ssec: FCM}, these length scales are truncated for points near the boundaries to account for the truncation of the nonlocal horizon. We have illustrated this effect on the length scales $l_A$ and $l_B$ in Fig.~(\ref{fig: beam_horizon}), where we have considered three points (P,Q, and R) at different distances from the beam boundaries. At points P and R close to the boundaries, the asymmetry in the nonlocal horizons is addressed by an appropriate truncation of the relevant length scales ($l_A^P$ and $l_B^R$). At point Q, which is sufficiently within the domain, a full horizon is captured by the complete length scales on either side of the point. Note that equal values for $l_A^Q$ and $l_B^Q$ represent an isotropic horizon.

The axial stress $\tilde{\sigma}_{11}$ resulting from the fractional-order strain is evaluated using Eq.~\eqref{eq: constt_isot} as:
\begin{equation}
    \label{eq: axial_stress_ebt}
    \tilde{\sigma}_{11}(x_1,x_3)=E\tilde{\epsilon}_{11}(x_1,x_3)
\end{equation}
where $E$ denotes the Young's Modulus of the isotropic elastic beam.
The deformation energy $\mathcal{U}$ accumulated in the nonlocal elastic beam due to the above stress-strain distribution is obtained as:
\begin{equation}
    \label{eq: def_energy}
    \mathcal{U}=\frac{1}{2}\int_{\Omega} \tilde{\sigma}_{11}(x_1,x_3)\tilde{\epsilon}_{11}(x_1,x_3)\mathrm{d}V
\end{equation}
where $\Omega$ denotes the volume occupied by the beam. Note that, by using the definition for the nonlocal strain in Eq.~(\ref{eq: trans_e13}) in the Euler-Bernoulli displacement field given in Eq.~(\ref{eq: euler_bern_beam}), a non-zero transverse shear strain would be obtained. The remaining strains are trivially zero for the chosen displacement field, analogous to classical formulations. The non-zero transverse shear strains are not considered in the deformation energy in Eq.~(\ref{eq: def_energy}) because for under slender beam assumptions, the rigidity to transverse shear forces is much higher when compared to the bending rigidity. Hence, the contribution of the transverse shear deformation towards the deformation energy of the solid can be neglected. 

The total potential energy functional, assuming that the beam is subjected to distributed forces $F_a(x_1)$ and $F_t(x_1)$ acting on the mid-plane along the axial and transverse directions, is given by:
\begin{equation}
    \label{eq: pot_energy}
    \Pi[\textbf{u}(\textbf{x})]=\mathcal{U}-\int_{0}^{L}F_a(x_1) u_0(x_1)\mathrm{d}x_1-\int_{0}^{L}F_t(x_1) w_0(x_1)\mathrm{d}x_1
\end{equation}
To obtain the total potential energy of the system, we have assumed no body forces to be applied. Using the above constitutive model for the fractional-order nonlocal and nonlinear elastic beam, we can derive the governing differential equations in their strong form and the associated boundary conditions by imposing optimality conditions on the total potential energy functional $\Pi[\textbf{u}(\textbf{x})]$.

\subsection{Governing equations}
Using a variational approach, the governing differential equations and the associated boundary conditions for the fractional-order beam are:
\begin{subequations}
\label{eq: all_governing_equations}
\begin{equation}\label{eq: axial_gde}
        {}^{R-RL}_{x_1-l_B}D^{\alpha}_{x_1+l_A}\left[N(x_1)\right]+F_a(x_1)=0~~\forall~~x_1\in(0,L)
\end{equation}
\begin{equation}\label{eq: transverse_gde}
        \frac{\mathrm{d}}{\mathrm{d}x_1}\big[{}^{R-RL}_{x_1-l_B}D^{\alpha}_{x_1+l_A}\left[M(x_1)\right]\big]+ {}^{R-RL}_{x_1-l_B}D^{\alpha}_{x_1+l_A} \left[N(x_1)D_{x_1}^{\alpha}\left[w_0(x_1)\right]\right]+F_t(x_1)=0~~\forall~~x_1\in(0,L)
\end{equation}
and subject to the boundary conditions:
\begin{equation}\label{eq: axial_bcs}
    N(x_1)=0~~\text{or}~~\delta u_0(x_1)=0~~\text{at}~~x_1\in\{0,L\}
\end{equation}
\begin{equation}\label{eq: transverse_moment_bcs}
    M(x_1)=0~~\text{or}~~\delta\Bigg[\frac{\mathrm{d}w_0(x_1)}{\mathrm{d}x_1}\Bigg]=0~~\text{at}~~x_1\in\{0,L\}
\end{equation}
\begin{equation}\label{eq: transverse_force_bcs}
    \frac{\mathrm{d}M(x_1)}{\mathrm{d}x_1}+N(x_1)\frac{\mathrm{d}w_0(x_1)}{\mathrm{d}x_1}=0~~\text{or}~~\delta w_0(x_1)=0~~\text{at}~~x_1\in\{0,L\}
\end{equation}
\end{subequations}
In the above, $N(x_1)$ and $M(x_1)$ are the axial and bending stress resultants, and are given by:
\begin{subequations}
\label{eq: stress_resultants}
\begin{equation}
\label{eq: stress_resultant_def_1}
    N(x_1)=\int_{-b/2}^{b/2}\int_{-h/2}^{h/2}\tilde{\sigma}_{11}~\mathrm{d}x_3~\mathrm{d}x_2
\end{equation}
\begin{equation}
\label{eq: stress_resultant_def_2}
    M(x_1)=\int_{-b/2}^{b/2}\int_{-h/2}^{h/2}x_3~\tilde{\sigma}_{11}~\mathrm{d}x_3~\mathrm{d}x_2
\end{equation}
\end{subequations}
where $b$ denotes the width of the beam. In the fractional-order governing equations provided in Eq.~\eqref{eq: all_governing_equations}, ${}^{R-RL}_{x_1-l_B}D^{\alpha}_{x_1+l_A}[\cdot]$ is a Riesz-type Riemann-Liouville (R-RL) derivative of order $\alpha\in(0,1)$. The R-RL fractional derivative of an arbitrary function $f(x_1)$ is defined analogous to the RC derivative given in Eq.~\eqref{eq: frac_der_def} as:
\begin{equation}
    \label{eq: r_rl_frac_der_def}
    {}^{R-RL}_{x_1-l_B}D^{\alpha}_{x_1+l_A}f(x_1)=\frac{1}{2}\Gamma(2-\alpha)\left[l_B^{\alpha-1}\left({}^{RL}_{x_1-l_B}D^{\alpha}_{x_1}f(x_1)\right)-l_A^{\alpha-1}\left({}^{RL}_{x_1}D^{\alpha}_{x_1+l_A}f(x_1)\right)\right]
\end{equation}
where ${}^{RL}_{x_1-l_B}D^{\alpha}_{x_1}f(x_1)$ and ${}^{RL}_{x_1}D^{\alpha}_{x_1+l_A}f(x_1)$ are the left- and right-handed Riemann Liouville (RL) derivatives of an arbitrary function $f(x_1)$ to the order $\alpha$, respectively. 
The detailed steps to obtain the strong form of the governing equations by minimization of the total potential energy are reported in Appendix 2.

Note that the governing equations for the axial and transverse displacements are coupled, similar to what is seen in the classical nonlinear Euler-Bernoulli beam formulation. This is unlike the linear elastic fractional-order response of a nonlocal beam, where the governing equations for the axial and transverse displacements are uncoupled \cite{patnaik2019FEM}. Further, as expected, the classical nonlinear Euler-Bernoulli beam governing equations and boundary conditions are recovered for $\alpha=1$.
As mentioned previously, the above form of the governing equations for the fractional-order beam do not admit closed-form analytical solutions for the most general loading and boundary conditions. Therefore, in the following, we develop a fractional-order finite element method (f-FEM) to obtain the numerical solution of the nonlinear fractional-order governing equations. We will then use the f-FEM to investigate the characteristic response of a nonlinear fractional-order beam.

\section{Fractional Finite Element Method (f-FEM)}
\label{sec: ffem}
The fractional-order nonlocal and nonlinear boundary value problem described by Eqs.~(\ref{eq: axial_gde}-\ref{eq: transverse_force_bcs}) are numerically solved using finite element method. The fractional-order boundary value problem builds upon the FE methods developed in the literature for integral models of nonlocal elasticity \cite{pisano2009nonlocal,norouzzadeh2017finite}. However, several modifications are necessary due to both the specific form of the attenuation functions used in the definition of the fractional-order derivatives, and to the nonlocal continuum model adopted in this study. Note that, although a f-FEM formulation was developed in \cite{patnaik2019FEM}, the study dealt exclusively with linear fractional-order boundary value problems, while here we present the approach in the context of a nonlinear BVP. We also highlight that, although the f-FEM is presented here for nonlinear fractional-order beams, its formulation can be easily extended to model the nonlocal and nonlinear behavior of plates and shells. 

\subsection{f-FEM formulation}
\label{sec: ffem_method}
We formulate the f-FEM starting from a discretized form of the total potential energy functional $\Pi[\textbf{u}(\textbf{x})]$ given in Eq.~(\ref{eq: pot_energy}).
For this purpose, the beam domain $\Omega=~[0,L]$ along the $x_1$ direction, is uniformly discretized into disjoint two-noded elements $\Omega^e_i=(x_1^i,x_1^{i+1})$, such that $\cup_{i=1}^{N_e}\Omega^e_i = \Omega$, and $\Omega^e_j\cap\Omega^e_k=\emptyset~\forall~j\neq k$. $N_e$ is the total number of discretized elements and the length of each element $\Omega^e\in \Omega$ is $l_e$.
The unknown displacement field variables $u_0(x_1)$ and $w_0(x_1)$ at any point $x_1 \in \Omega^e_i$ are evaluated by interpolating the corresponding nodal degrees of freedom for $\Omega^e_i$. From the definitions given in Eq. \eqref{eq: von_karman_ebt}, it is clear that the interpolation of the axial and transverse displacement fields following Euler-Bernoulli beam theory would require at least Lagrangian ($\mathcal{L}_i$, $i=1,2$) and Hermitian ($\mathcal{H}_j$, $j=1,2,3,4$) interpolation functions, respectively. Consequently, the axial $u_0(x_1)$ and transverse $w_0(x_1)$ displacement fields can be written  as:
\begin{subequations}
\label{eq: combined_interpolations}
\begin{equation}
    \label{eq: intpl_disp_u}
    \{u_0(x_1)\}=[S_u(x_1)]\{U_l (x_1)\}
\end{equation}
\begin{equation}
    \label{eq: intpl_disp_w}
    \{w_0(x_1)\}=[S_w(x_1)]\{W_l (x_1)\}
\end{equation}
where the localized displacement coordinate vectors $\{U_l (x_1)\}$ and $\{W_l (x_1)\}$ corresponding to the two-noded 1D element enclosing the point $x_1$ are given as:
\begin{equation}
    \{U_l(x_1)\}^{T}=\begin{bmatrix}
        u_0^{(1)}&u_0^{(2)}
        \end{bmatrix}
\end{equation}
\begin{equation}
    \{W_l(x_1)\}^{T}=\begin{bmatrix}
        w_0^{(1)}&\frac{\mathrm{d}w_0^{(1)}}{\mathrm{d}x_1}&w_0^{(2)}&\frac{\mathrm{d}w_0^{(2)}}{\mathrm{d}x_1}
        \end{bmatrix}
\end{equation}
and the associated shape function matrices are:
\begin{equation}
    [S_u(x_1)]=\begin{bmatrix}
    \mathcal{L}_1^e(x_1) & \mathcal{L}_2^e(x_1)
    \end{bmatrix}
\end{equation}
\begin{equation}
    [S_w(x_1)]=\begin{bmatrix}
     \mathcal{H}_1^e(x_1) & \mathcal{H}_2^e(x_1) & \mathcal{H}_3^e(x_1) & \mathcal{H}_4^e(x_1)
    \end{bmatrix}
\end{equation}
\end{subequations}
The superscripts $(\cdot)^{(1)}$ and $(\cdot)^{(2)}$ in the above equation denote the local node numbers for the two-noded element. Using the definition of the RC derivative in Eq.~(\ref{eq: frac_der_def}), the fractional derivative $D_{x_1}^{\alpha} \left[u_0(x_1)\right]$ is expressed as:
\begin{equation}
    \label{eq: rc_disp_simpl}
    D_{x_1}^{\alpha} \left[u_0(x_1)\right]=\frac{1}{2}(1-\alpha)\left[l_A^{\alpha-1}\int_{x_1-l_A}^{x_1}\frac{D^1[u_0(s_1)]}{(x_1-s_1)^{\alpha}}~\mathrm{d}s_1+l_B^{\alpha-1}\int_{x_1}^{x_1+l_B} \frac{D^1[u_0(s_1)]}{(s_1-x_1)^{\alpha}}~ \mathrm{d}s_1\right]
\end{equation}
where $D^1(\cdot)$ is the first integer-order derivative and $s_1$ is a dummy variable along the direction $\textbf{x}_1$ used within the definition of the fractional-order derivative. 
The above expression can be recast as:
\begin{subequations}
\label{eq: rc_simpl_form_combined}
\begin{equation}
    \label{eq: rc_simpl_form2}
    D_{x_1}^{\alpha} \left[u_0(x_1)\right]=\int_{x_1-l_A}^{x_1} A(x_1,s_1,l_A,\alpha)~D^1 \left[u_0(s_1)\right]~\mathrm{d}s_1+\int_{x_1}^{x_1+l_B} A(x_1,s_1,l_B,\alpha)~D^1\left[u_0(s_1)\right]~\mathrm{d}s_1
\end{equation}
where
\begin{equation}
    \label{eq: frac_kernel_gen}
    A(x,y,l,\alpha)=\frac{(1-\alpha)}{2}~l^{\alpha-1}~\frac{1}{|x-y|^{\alpha}}
\end{equation}
\end{subequations}
is the kernel of the fractional derivative. 
Note that kernel $A(x,y,l,\alpha)$ is a function of the relative distance between the points $x$ and $y$, and can be interpreted as an analog form of the attenuation functions typically used in integral models of nonlocal elasticity. Clearly, the attenuation decays as a power-law in the distance with an exponent equal to the order $\alpha$ of the fractional derivative. Similar expressions for $D^{\alpha}_{x_1}\left[w_0(x_1)\right]$ and $D^{\alpha}_{x_1}\left[\mathrm{d}w_0(x_1)/\mathrm{d}x_1\right]$ can also be derived using an equivalent approach.

Note that $D^{\alpha}_{x_1}\left[u_0(x_1)\right]$, $D^{\alpha}_{x_1}\left[w_0(x_1)\right]$ and $D^{\alpha}_{x_1}\left[\mathrm{d}w_0(x_1)/\mathrm{d}x_1\right]$ contain the integer-order derivatives $D^1_{s_1}[u_0(s_1)]$, $D^1_{s_1}[w_0(s_1)]$ and $D^1_{s_1}[\mathrm{d}w_0(s_1)/\mathrm{d}s_1]$ (equivalent to $D^2_{s_1} [w_0(s_1)]$). These integer-order derivatives at $s_1$ are evaluated in terms of the nodal variables corresponding to the element $\Omega_p^e$, such that $s_1\in \Omega_p^e$. Using Eq.~\eqref{eq: combined_interpolations}, the integer order derivatives can be expressed as:
\begin{equation}
\label{eq: deriv_disp_rel}
    D^1[u_0(s_1)]=[B_u(s_1)]\{U_l(s_1)\}, ~~~D^1[w_0(s_1)]=[B_w(s_1)]\{W_l(s_1)\},
    ~~~D^2[w_0(s_1)]=[B_\theta(s_1)]\{W_l(s_1)\}
\end{equation}
where the $[B_{\square}]$ ($\square\in\{u,w,\theta\}$) matrices are expressed in terms of the associated shape functions for the two-noded element as:
\begin{subequations}
\label{eq: b_mats}
\begin{equation}
    [B_u(s_1)]=\begin{bmatrix}
    \frac{\mathrm{d}\mathcal{L}^e_1}{\mathrm{d}s_1} & \frac{\mathrm{d}\mathcal{L}^e_2}{\mathrm{d}s_1}
    \end{bmatrix}
    \end{equation}
\begin{equation}
    [B_w(s_1)]=\begin{bmatrix} \frac{\mathrm{d}\mathcal{H}_1^e}{\mathrm{d}s_1} & \frac{\mathrm{d}\mathcal{H}_2^e}{\mathrm{d}s_1} & \frac{\mathrm{d}\mathcal{H}_3^e}{\mathrm{d}s_1} & \frac{\mathrm{d}\mathcal{H}_4^e}{\mathrm{d}s_1}
    \end{bmatrix}
\end{equation}
\begin{equation}
    [B_\theta(s_1)]=\begin{bmatrix}
     \frac{\mathrm{d}^2\mathcal{H}_1^e}{\mathrm{d}s_1^2} & \frac{\mathrm{d}^2\mathcal{H}_2^e}{\mathrm{d}^2 s_1} &  \frac{\mathrm{d}^2\mathcal{H}_3^e}{\mathrm{d}^2 s_1} & \frac{\mathrm{d}^2\mathcal{H}_4^e}{\mathrm{d}^2 s_1}
    \end{bmatrix}
\end{equation}
\end{subequations}
Using the above expression for the integer-order derivative $D^{1}[u_0(x_1)]$, the fractional-order derivative in Eq.~\eqref{eq: rc_simpl_form_combined} can be written as:
\begin{subequations}
\begin{equation}
    \label{eq: rc_disp_simpl_3}
    D_{x_1}^{\alpha} \left[u_0(x_1)\right]=\int_{x_1-l_A}^{x_1+l_B} A_p(x_1,s_1,l_A,l_B,\alpha)[B_u(s_1)]\{U_l(s_1)\}\mathrm{d}s_1
\end{equation}
where the following definition for the attenuation function $A_p(x_1,s_1,l_A,l_B,\alpha)$ is employed:
\begin{equation}
\label{eq: atten_func}
    A_p(x_1,s_1,l_A,l_B,\alpha)=\begin{cases}
    A(x_1,s_1,l_A,\alpha)& ~~ s_1\in{(x_1-l_A,x_1)},\\
A(x_1,s_1,l_B,\alpha) & ~~ s_1\in{(x_1,x_1+l_B)}.
    \end{cases}
\end{equation}
\end{subequations}
Similar expressions are derived for $D^\alpha_{x_1} [w_0(x_1)]$ and $D^\alpha_{x_1} [\mathrm{d}w_0(x_1)/\mathrm{d}x_1]$:
\begin{subequations}
\begin{equation}
    D_{x_1}^{\alpha} \left[w_0(x_1)\right]=\int_{x_1-l_A}^{x_1+l_B} A_p(x_1,s_1,l_A,l_B,\alpha)[B_w(s_1)]\{W_l(s_1)\}\mathrm{d}s_1
\end{equation}
\begin{equation}
    D_{x_1}^{\alpha} \left[\frac{\mathrm{d}w_0(x_1)}{\mathrm{d}x_1}\right]=\int_{x_1-l_A}^{x_1+l_B} A_p(x_1,s_1,l_A,l_B,\alpha)[B_\theta(s_1)]\{W_l(s_1)\}\mathrm{d}s_1
\end{equation}
\end{subequations}
Note that the evaluation of the fractional derivative requires a convolution of the integer-order derivatives across the horizon of nonlocality $(x_1-l_A,x_1+l_B)$. 
Thus, while obtaining the FE approximation in Eq.~(\ref{eq: rc_disp_simpl_3}), the nonlocal contributions from the different finite elements in the horizon must be properly attributed to the corresponding nodes. In order perform this mapping of the nonlocal contributions from elements in the horizon, we transform the nodal values $\{U_l(s_1)\}$ and $\{W_l(s_1)\}$ into the respective global variable vectors $\{U_g\}$ and $\{W_g\}$ with the help of connectivity matrices in the following manner:
\begin{subequations}
\label{eq: conversion_to_global_form}
\begin{equation}
    \{U_l(s_1)\}=[\tilde{\mathcal{C}_u}(x_1,s_1)]\{U_g\}
\end{equation}
\begin{equation}
    \{W_l(s_1)\}=[\tilde{\mathcal{C}_w}(x_1,s_1)]\{W_g\}
\end{equation}
\end{subequations}
The connectivity matrices $[\tilde{\mathcal{C}_u}(x_1,s_1)]$ and  $[\tilde{\mathcal{C}_w}(x_1,s_1)]$ are defined for the axial and transverse displacement vectors, respectively. These matrices are designed such that they are non-zero only if $s_1$ lies in the domain $(x_1-l_A,x_1+l_B)$, the horizon of nonlocality for $x_1$. It is evident that these matrices activate the contribution of the nodes enclosing $s_1$ for the numerical evaluation of the convolution integral given in Eq. \eqref{eq: rc_disp_simpl_3} \cite{patnaik2019FEM}. Using the above formalism, Eq. \eqref{eq: rc_disp_simpl_3} can be recast as:
\begin{subequations}
\label{eq: frac_der_axial}
\begin{equation}
\label{eq: frac_der_u_final}
    D^{\alpha}_{x_1}\left[u_0(x_1)\right]=[\tilde{B}_{u}(x_1)]\{U_g\}
\end{equation}
where
\begin{equation}
\label{eq: frac_der_Bu_final}
    [\tilde{B}_u(x_1)]=\int_{x_1-l_A}^{x_1+l_B}A_p(x_1,s_1,l_A,l_B,\alpha)[B_u(s_1)][\tilde{\mathcal{C}_u}(x_1,s_1)]~\mathrm{d}s_1
\end{equation}
\end{subequations}
Similarly, we express $D^{\alpha}_{x_1}\left[w_0(x_1)\right]$ and $ D^{\alpha}_{x_1}\left[{\mathrm{d}w_0(x_1)}/{\mathrm{d}x_1}\right]$ as:
\begin{subequations}
\label{eq: frac_der_trans}
\begin{equation}
\label{eq: frac_der_w_final}
    D^{\alpha}_{x_1}\left[w_0(x_1)\right]=[\tilde{B}_{w}(x_1)]\{W_g\}
\end{equation}
\begin{equation}
\label{eq: frac_der_theta_final}
    D^{\alpha}_{x_1}\left[\frac{\mathrm{d}w_0(x_1)}{\mathrm{d}x_1}\right]=[\tilde{B}_{\theta}(x_1)]\{W_g\}
\end{equation}
where
\begin{equation}
\label{eq: frac_der_Bw_final}
    [\tilde{B}_w(x_1)]=\int_{x_1-l_A}^{x_1+l_B}A_p(x_1,s_1,l_A,l_B,\alpha)[B_w(s_1)][\tilde{\mathcal{C}_w}(x_1,s_1)]~\mathrm{d}s_1
\end{equation}
\begin{equation}
\label{eq: frac_der_Btheta_final}
    [\tilde{B}_\theta(x_1)]=\int_{x_1-l_A}^{x_1+l_B}A_p(x_1,s_1,l_A,l_B,\alpha)[B_\theta(s_1)][\tilde{\mathcal{C}_w}(x_1,s_1)]~\mathrm{d}s_1
\end{equation}
\end{subequations}
Note that $[B_u(s_1)]$, $[B_w(s_1)]$ and $[B_\theta(s_1)]$ have been given in Eq.~\eqref{eq: b_mats}.
\subsection{Nonlinear f-FEM model}
In the following, we will use the above formalism to obtain the algebraic equations of motion for the geometrically nonlinear response of the fractional-order nonlocal beam by minimizing the total potential energy of the beam. The first variation of the deformation energy $\mathcal{U}$ defined in Eq.~\eqref{eq: def_energy} is given by:
\begin{equation}
\label{eq: def_energy_beam}
    \delta \mathcal{U}=b\int_{0}^{L} \int_{-h/2}^{h/2} \delta \tilde{\epsilon}_{11}(x_1,x_3)~ \tilde{\sigma}_{11} (x_1,x_3) \mathrm{d}x_3 \mathrm{d}x_1
\end{equation}
Using the above expression for $\delta \mathcal{U}$ along with the nonlocal axial strain derived in Eq. \eqref{eq: von_karman_ebt} for the Euler-Bernoulli beam theory and the definitions of the stress resultants given in Eq. \eqref{eq: stress_resultants}, we obtain the following expression for the statement of the minimum potential energy $\delta \Pi =0$:
\begin{equation}
    \label{eq: pot_energy_simp}
    \begin{split}
    \int_{0}^{L} \left\{ \left(D_{x_1}^{\alpha}\left[\delta u_0(x_1)\right]+D_{x_1}^{\alpha} \left[w_0(x_1)\right]~D_{x_1}^{\alpha}\left[\delta w_0(x_1)\right]\right)N(x_1)-D_{x_1}^{\alpha}\left[\frac{\mathrm{d}\delta w_0(x_1)}{\mathrm{d}x_1}\right]M(x_1)\right\}\mathrm{d}x_1&\\
    -\int_0^L F_t(x_1)\delta w_0(x_1) \mathrm{d}x_1-\int_0^L F_a(x_1)\delta u_0(x_1) \mathrm{d}x_1&=0
    \end{split}
\end{equation}
Note that the variations $\delta u_0(x_1)$ and $\delta w_0(x_1)$ are independent of each other. Now, by collecting the terms containing $\delta u_0(x_1)$ and $\delta w_0(x_1)$, the above weak form is expressed in the form of the following equations:
\begin{subequations}
\label{eq: weak_form_eq}
\begin{equation}
    \label{eq: weak_form_axial_eq}
    \int_{0}^{L} \left\{D_{x_1}^{\alpha} \left[\delta u_0(x_1)\right] N(x_1)-F_a(x_1)\delta u_0(x_1)\right\}\mathrm{d}x_1=0
\end{equation}
\begin{equation}
    \label{eq: weak_form_trans_eq}
    \int_0^L \left\{\left(D_{x_1}^{\alpha} \left[w_0(x_1)\right]~D_{x_1}^{\alpha}\left[\delta w_0(x_1)\right]\right)N(x_1)-D_{x_1}^{\alpha}\left[\frac{\mathrm{d}\delta w_0(x_1)}{\mathrm{d}x_1}\right]M(x_1)-F_t(x_1)\delta w_0(x_1)\right\}\mathrm{d}x_1=0
\end{equation}
\end{subequations}
We recast the above equations in terms of of the displacement field variables by using the expressions for the stress resultants given in Eqs.~(\ref{eq: constt_isot}) and (\ref{eq: stress_resultants}), and obtain the following:
\begin{subequations}
\label{eq: weak_form_eq_1}
\begin{equation}
    \label{eq: weak_form_eq_axial_1}
    \int_{0}^{L}\Bigg\{A_{11}\left(D_{x_1}^{\alpha}\left[\delta u_0(x_1)\right]\right) \big(D_{x_1}^{\alpha} \left[u_0(x_1)\right]+\frac{1}{2}\left(D_{x_1}^{\alpha} \left[w_0(x_1)\right]\right)^2\big)-F_a(x_1)\delta u_0(x_1)\Bigg\} \mathrm{d}x_1=0
\end{equation}
\begin{equation}
\label{eq: weak_form_eq_trans_1}
\begin{split}
    \int_{0}^{L} &\Bigg\{ A_{11}~\bigg(D_{x_1}^{\alpha} \left[w_0(x_1)\right]~D_{x_1}^{\alpha}\left[\delta w_0(x_1)\right]\bigg)~ \big(D_{x_1}^{\alpha} \left[u_0(x_1)\right]+\frac{1}{2}\left(D_{x_1}^{\alpha} \left[w_0(x_1)\right]\right)^2\big)\\
    &+D_{11}~ \left(D_{x_1}^{\alpha}\left[\frac{\mathrm{d}\delta w_0(x_1)}{\mathrm{d}x_1}\right]\right)~\left(D_{x_1}^{\alpha}\left[\frac{\mathrm{d}w_0(x_1)}{\mathrm{d}x_1}\right]\right)-F_t(x_1)~\delta w_0(x_1) \Bigg\}~\mathrm{d}x_1=0
\end{split}
\end{equation}
\end{subequations}
where $A_{11}=Ebh$ and $D_{11}=Ebh^3/12$ are the axial and bending stiffness of the beam, respectively.
The above equations in the weak form are converted into algebraic equations using the f-FEM model discussed in \S\ref{sec: ffem_method}. Now, by using Eqs.~(\ref{eq: frac_der_axial}-\ref{eq: frac_der_trans}) and enforcing minimization, the above integrals reduce to the following algebraic equations:
\begin{equation}
    \label{eq: algebraic_eqs}
    \begin{bmatrix} [\tilde{K}_{11}] & [\tilde{K}_{12}]\\ [\tilde{K}_{21}] & [\tilde{K}_{22}] \end{bmatrix}
    \begin{Bmatrix} \{U_g\}\\ \{W_g\} \end{Bmatrix}=
    \begin{Bmatrix} \{F_A\}\\ \{F_T\} \end{Bmatrix}
\end{equation}
In the above system of equations, the different stiffness matrices are obtained as:
\begin{subequations}
    \label{eq: stiffness_mats}
\begin{equation}
    \label{eq: stiffness_mats_k11}
    [\tilde{K}_{11}]=\int_{0}^{L}A_{11} [\tilde{B}_{u}(x_1)]^{T} [\tilde{B}_{u}(x_1)]~\mathrm{d}x_1
\end{equation}
\begin{equation}
    \label{eq: stiffness_mats_k12}
    [\tilde{K}_{12}]=\frac{1}{2}\int_{0}^{L}A_{11} \left(D_{x_1}^{\alpha}[w_0(x_1)]\right) [\tilde{B}_{u}(x_1)]^{T} [\tilde{B}_{w}(x_1)]~\mathrm{d}x_1
\end{equation}
\begin{equation}
    \label{eq: stiffness_mats_k21}
    [\tilde{K}_{21}]=\int_{0}^{L}A_{11} \left(D_{x_1}^{\alpha}[w_0(x_1)]\right) [\tilde{B}_{w}(x_1)]^{T} [\tilde{B}_{u}(x_1)]~\mathrm{d}x_1
\end{equation}
\begin{equation}
    \label{eq: stiffness_mats_k22}
    [\tilde{K}_{22}]=\int_{0}^{L}D_{11} [\tilde{B}_{\theta}(x_1)]^{T} [\tilde{B}_{\theta}(x_1)]~\mathrm{d}x_1+\frac{1}{2}\int_0^L  \left[A_{11}\left(D_{x_1}^{\alpha}[w_0(x_1)]\right)^2\right][\tilde{B}_w(x_1)]^T[\tilde{B}_w(x_1)]~\mathrm{d}x_1
\end{equation}
and the axial and transverse nodal force vectors are obtained as:
\begin{equation}
    \{F_A\}^{T}=\int_{0}^{L}F_a(x_1)~[S_u(x_1)]~\mathrm{d}x_1
\end{equation}
\begin{equation}
    \{F_T\}^{T}=\int_{0}^{L}F_t(x_1)~[S_w(x_1)]~\mathrm{d}x_1
\end{equation}
\end{subequations}
The solution of the system of algebraic equations in Eq.~\eqref{eq: algebraic_eqs} provides the nodal generalized displacement coordinates, which can then be used along with the Euler-Bernoulli relations in Eq.~\eqref{eq: euler_bern_beam} to determine the displacement field at any point of the beam. The effect of the geometric nonlinearity is evident from the expressions of the stiffness matrices given in Eq. \eqref{eq: stiffness_mats}, where it appears that the stiffness matrices are a function of the deformed system configuration.
Note that the linear f-FEM model can be obtained by ignoring the contribution of $D_{x_1}^{\alpha}[w_0(x_1)]$ in Eq. \eqref{eq: stiffness_mats}. Further, similar to classical beam models accounting for geometric nonlinearity, a coupling between the axial and bending responses is obtained through the non-zero $[\tilde{K}_{12}]$ and $[\tilde{K}_{21}]$ matrices.

In this study, we solve the nonlinear algebraic equations using an incremental-iterative Newton-Raphson (NR) numerical scheme. Similar to classical nonlinear models, the NR procedure for the fractional-order nonlinear equations also requires the evaluation of the tangent stiffness matrix. The procedure to derive the tangent stiffness matrix $[\tilde{K}_T^i]$ corresponding to the $i$-th iteration of the NR method follows from the classical approach:
\begin{equation}
    \label{eq: tangent_stiffness_matrix}
    [\tilde{K}_T^{i}]=\left[\frac{\partial \{R\}}{\partial \{X\}}\right]^{i}
\end{equation}
where $\{X\}$ is a vector containing all the displacement coordinates of the beam and $\{R\}$ is the residual vector evaluated at the current iteration. The residual $\{R\}$ at the $i$-th iteration step is given by: 
\begin{equation}
    \{R(\{X\})\}^{i}=[\tilde{K}_S^i]\{X\}^i-\{F\}
\end{equation}
where $\tilde{K}_S^i$ denotes the stiffness matrix of the entire system. In order to obtain the solution of the nonlinear equations, we also adopt a load increment procedure wherein the nonlocal response of the beam at each load step is evaluated using the NR iteration:
\begin{equation}
    \label{eq: NR_procedure}
    \{X\}^{i+1}=\{X\}^{i}-\left[\tilde{K}_T^{i}\right]^{-1}\{R(\{X\})\}^{i}
\end{equation} The iterations at each load level are continued until the residual becomes less than a chosen tolerance.
\subsection{Numerical integration for nonlocal matrices}
In this following, we provide the details of the numerical scheme used to integrate both the stiffness and the tangent stiffness matrices. Note that the numerical procedure to evaluate the force vector follows directly from classical FE formulations. Hence, we do not provide a complete account of all the steps.
Evaluation of the stiffness matrices for the nonlocal system given in Eq. \eqref{eq: stiffness_mats} requires the evaluation of the nonlocal matrices $[\tilde{B}_{\square}]$ $(\square \in \{u,w,\theta\})$. As seen in Eqs.~(\ref{eq: frac_der_axial}-\ref{eq: frac_der_trans}), this involves a convolution of the integer-order derivatives with the fractional-order attenuation function $A_{p}(x_1,s_1,l_A,l_B,\alpha)$ over the nonlocal horizon $(x_1-l_A,x_1+l_B)$. Clearly, the numerical model for fractional-order derivative involves additional integration over the nonlocal horizon to account for nonlocal interactions.
These nonlocal interactions across the horizon of nonlocality have already been accounted for numerically in \cite{polizzotto2001nonlocal,pisano2009nonlocal}. However, differently from these studies, the attenuation function $A_{p}(x_1,s_1,l_A,l_B,\alpha)$ in the fractional-order model involves an end-point singularity (more specifically, at $x_1=s_1$) due to the nature of the kernel (see Eq.~(\ref{eq: atten_func})). The fractional-order nonlocal interactions as well as the end-point singularity were addressed in \cite{patnaik2019FEM} and are briefly reviewed here. 

In the following, we describe the procedure for the numerical evaluation of the nonlinear stiffness matrix $[\tilde{K}_{12}]$. The same procedure can be extended to evaluate the other stiffness matrices, hence their discussion is not provided here. In order to perform the numerical integration of the nonlocal and nonlinear stiffness matrix $[\tilde{K}_{12}]$, we adopt an isoparametric formulation and introduce a natural coordinate system $\xi_1$. The Jacobian of the transformation $x_1\rightarrow\xi_1$ is given as $J(\xi_1)$. Now, by using the Gauss-Legendre quadrature rule, the matrix $[\tilde{K}_{12}]$ is approximated as:
\begin{equation}
    [\tilde{K}_{12}]
    \approx \sum\limits_{i=1}^{N_e}\sum\limits_{j=1}^{N_{GP}} \hat{w}_j J^i\left\{A_{11}\left(D_{x_1}^{\alpha}\left[w_0(\xi_1^{i,j})\right]\right)\right\}[\tilde{B}_u(\xi_1^{i,j})]^T[\tilde{B}_w(\xi_1^{i,j})]
\end{equation}
where $\xi_1^{i,j}$ is the $j-$th Gauss integration point in the $i-$th element, $\hat{w}_j$ is the corresponding weight for numerical integration, and $N_{GP}$ is the total number of Gauss points chosen for the numerical integration, such that $j\in \{1,...N_{GP}\}$. $J^i$ is the Jacobian of the coordinate transformation for the $i-$th element. The hat symbol on the weight is used to distinguish it from the transverse displacement.
As previously highlighted, the matrices $[\tilde{B}_{\square}]$ ($\square\in\{u,w,\theta\}$) (equivalently, $[\tilde{B}_{\square}(\xi_1^{i,j})]$ in the discretized form) involve an additional integration (see Eq.~(\ref{eq: frac_der_Bu_final})) due to the fractional-order nonlocality and are evaluated in the following manner:
\begin{equation}
    \left[\tilde{B}_\square(x_1^{i,j})\right]=\int_{x_1^{i,j}-l_A}^{x_1^{i,j}+l_B}A_p(x_1^{i,j},s_1,l_A,l_B,\alpha)[B_\square(s_1)][\tilde{\mathcal{C}_\square}(x_1^{i,j},s_1)]~\mathrm{d}s_1
\end{equation}
where $x_1^{i,j}$ is the Cartesian coordinate of the Gauss point $\xi_{1}^{i,j}$ and $[B_{\square}(s_1)]$ is given in Eq.~\eqref{eq: b_mats}. Using Eq.~\eqref{eq: atten_func}, we obtain the following expression for $\left[\tilde{B}_\square(x_1^{i,j})\right]$:
\begin{subequations}
\begin{equation}
\label{eq: b_mat_numer_simp1}
    \left[\tilde{B}_\square(x_1^{i,j})\right]=\int_{x_1^{i,j}-l_A}^{x_1^{i,j}}\mathcal{I}_L\mathrm{d}s_1+\int_{x_1^{i,j}}^{x_1^{i,j}+l_B}\mathcal{I}_R\mathrm{d}s_1
\end{equation}
where, the integrands $\mathcal{I}_L$ and $\mathcal{I}_R$ are given as:
\begin{equation}
    \mathcal{I}_L=A(x_1^{i,j},s_1,l_A,\alpha)[B_\square(s_1)][\tilde{\mathcal{C}_u}(x_1^{i,j},s_1)]
\end{equation}
\begin{equation}
    \mathcal{I}_R=A(x_1^{i,j},s_1,l_B,\alpha)[B_\square(s_1)][\tilde{\mathcal{C}_u}(x_1^{i,j},s_1)]
\end{equation}
\end{subequations}
Note that the limits of the integrals in Eq.~\eqref{eq: b_mat_numer_simp1} span over the elements that constitute the nonlocal horizon $(x_1^{i,j}-l_A,x_1^{i,j}+l_B)$ at $x_1^{i,j}$. These limits for the convolution integral follow from the terminals of the RC fractional derivative defined for the study of nonlocal continuum in Eq.~\eqref{eq: frac_der_def}. These integrals are evaluated numerically in the following manner:
\begin{subequations}
\label{eq: b_mat_int_Scheme}
\begin{equation}
    \int_{x_1^{i,j}-l_A}^{x_1^{i,j}} \mathcal{I}_L \mathrm{d}s_1\approx \underbrace{\int_{x_1^{i-N_A^{inf}}}^{x_1^{i-(N_A^{inf}-1)}} \mathcal{I}_L \mathrm{d}s_1+...\int_{x_1^{i-1}}^{x_1^{i}} \mathcal{I}_L \mathrm{d}s_1}_{\text{Gauss-Legendre Quadrature}}+ \underbrace{\int_{x_1^{i}}^{x_1^{i,j}} \mathcal{I}_L \mathrm{d}s_1}_{\text{Singularity at } x_1^{i,j}}
\end{equation}
\begin{equation}
    \int_{x_1^{i,j}}^{x_1^{i,j}+l_B} \mathcal{I}_R \mathrm{d}s_1\approx \underbrace{\int_{x_1^{i,j}}^{x_1^{i+1}} \mathcal{I}_R \mathrm{d}s_1}_{\text{Singularity at } x_1^{i,j}}+ \underbrace{\int_{x_1^{i+1}}^{x_1^{i+2}} \mathcal{I}_R \mathrm{d}s_1...+\int_{x_1^{i+(N_B^{inf}-1)}}^{x_1^{i+N_B^{inf}}} \mathcal{I}_R \mathrm{d}s_1}_{\text{Gauss-Legendre Quadrature}}
\end{equation}
\end{subequations}
In the above expressions, $N_A^{inf}$ and $N_B^{inf}$ are the number of (complete) elements in the nonlocal horizon to the left and right side of $x_1$, respectively. More specifically, $N_A^{inf}=\ceil{l_A/l_e}$ and ${N_B^{inf}}=\lfloor{l_B/l_e}\rfloor$ where $l_e$ is the length of the discretized element. The ceil ($\ceil{\cdot}$) and floor ($\lfloor{\cdot}\rfloor$) functions are used to round the number of elements to the greater integer on the left side and the lower integer on the right side. For $x_1^{i,j}$ close to the boundaries of the beam ($x_1=\{0,L\}$), $N_A^{inf}$ and $N_B^{inf}$ are truncated in order to account for asymmetric horizon lengths. This is essential to satisfy frame-invariance as discussed in \S\ref{sec: eb_model}.

\begin{figure}[h!]
    \centering
    \includegraphics[width=0.8\textwidth]{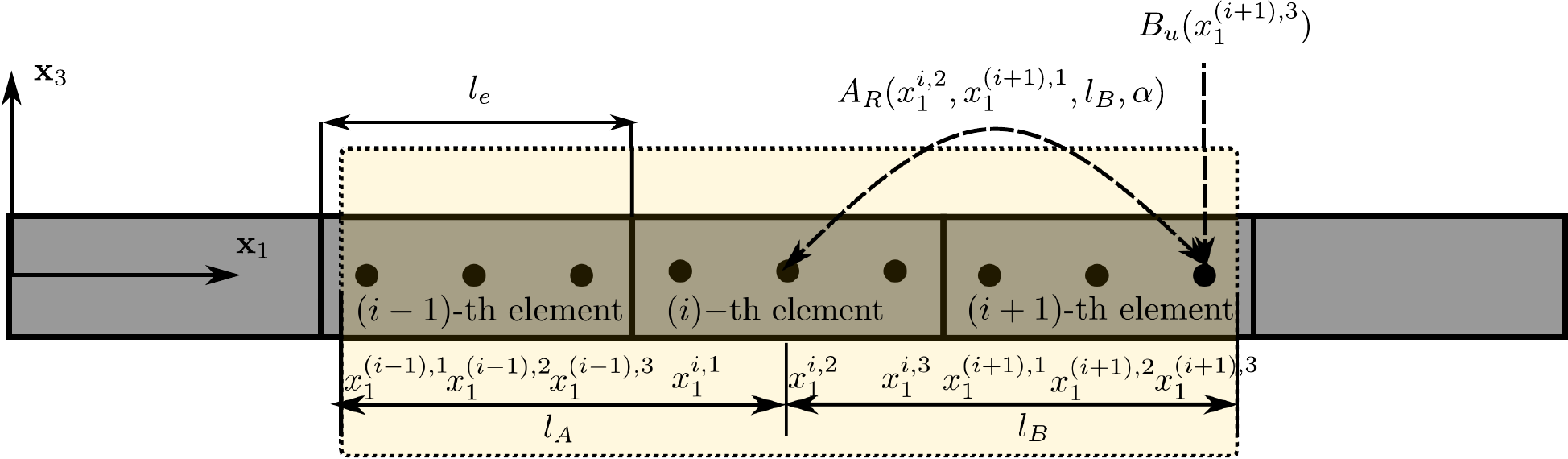}
    \caption{Illustration of the influence zone for the Gauss point $x_1^{i,2}$ and the evaluation of $[\tilde{B_u}(x_1^{i,2})]$.}
    \label{fig: nlfem_scheme}
\end{figure}

As discussed previously, due to the nature of the kernel of the fractional-order derivative, an end-point singularity occurs in the integrals at the Gauss point $x_1^{i,j}$ in the element $\Omega_i^e$. This is evident from the definitions of the left and right integrals given in Eq. \eqref{eq: b_mat_int_Scheme}. Following \cite{patnaik2019FEM}, this end-point singularity is circumvented by an analytical evaluation of these integrals over the elements containing the singularities. This analytical evaluation can be carried out by using the expressions for $[\tilde{B}_\square]$ ($\square\in u,w,\theta$) given in Eq.~(\ref{eq: frac_der_Bu_final}).
The integrals over the remaining elements (i.e. those without singularities) are evaluated using the Gauss-Legendre quadrature method. The expression for this quadrature based integration corresponding to the nonlocal contribution of the $r-$th element in the horizon of nonlocality of the Gauss point $x_1^{i,j}$ is given by:
\begin{equation}
\label{nonlocal_integral_example}
\begin{split}
        \int_{x_1^{r}}^{x_1^{r+1}}A_p(x^{i,j}_1,s_1,l_A,l_B,\alpha)[B_u(s_1)]&[\tilde{\mathcal{C}_u}(x_1^{i,j},s_1)]~\mathrm{d}s_1=\\
    &\sum_{k=1}^{N_{GP}}\hat{w}_k J^r A_p(x^{i,j}_1,{x}_1^{r,k},l_A,l_B,\alpha)[B_u({s}_1^{r,k})][\tilde{\mathcal{C}_u}(x_1^{i,j},{x}_1^{r,k})]
\end{split}
\end{equation}
where $x_1^{r,k}$ is the Cartesian coordinate of the $k-$th Gauss point in the $r-$th element along $\textbf{x}_1$, $\hat{w}_k$ is the corresponding weight, and $J^r$ is the Jacobian of the transformation for $r-$th element. The above numerical procedure is schematically illustrated in Fig. \ref{fig: nlfem_scheme}. We highlight that $A_p$ is a function of the Cartesian coordinates, therefore global coordinates should be used in its evaluation. 

As mentioned previously, the procedure followed for the evaluation of the $[\tilde{K}_{12}]$ stiffness matrix extends directly to the other stiffness matrices defined in Eq. \eqref{eq: stiffness_mats}, as well as the tangent stiffness matrix defined in Eq.~(\ref{eq: tangent_stiffness_matrix}). Note that all these matrices are pre-assembled in the global form due to the use of the connectivity matrices in Eq.~(\ref{eq: conversion_to_global_form}). The resulting matrices are then used in the NR iteration procedure outlined in Eq.~(\ref{eq: NR_procedure}) to obtain the generalized displacement coordinates for the geometrically nonlinear response of the fractional-order nonlocal beam.
\section{Application and Numerical Results}
\label{sec: num_results}
We use the f-FEM developed in \S\ref{sec: ffem} to study the geometrically nonlinear and fractional-order nonlocal elastic response of isotropic beams. More specifically, the effect of the fractional order $\alpha$ and of the length scales $l_A$, $l_B$ on the nonlocal elastic behavior of the nonlinear beam are studied. For all the numerical studies reported in this section, the aspect ratio of the beam is chosen to be $L/h=100$ in order to satisfy the slender beam assumption of the Euler-Bernoulli beam theory. The length of the beam is maintained at $L=1$ m, the width is kept constant at unity, and the elastic modulus is set at $E=3$ GPa. It is assumed here that that the length scales $l_A$ and $l_B$ for the isotropic beam are equal to a constant $l_f$ for points sufficiently within the domain of the beam. These length scales are truncated appropriately for points close to the beam boundaries as discussed in \S\ref{sec: eb_model} (see Fig.~(\ref{fig: beam_horizon})). Note that, although the presented f-FEM is capable of simulating the static response of the fractional-order nonlocal beam under both axial and transverse loads, here we  present only the results for transverse loads. Both the case of a uniformly distributed load per unit length (UDL) (in N/m) applied along the length of the beam, and of a point load (PtL) (in N) applied at a specified point. The transverse displacement, normal stress $\sigma_{11}$ and applied load are reported as:
\begin{equation}
    \overline{w}_0=\frac{w_0}{h};~~~~~\overline{\sigma}_{11}=\sigma_{11}\times\left(\frac{h}{L}\right)^2\frac{1}{q_0};~~~~~\overline{q}=q_0\times \frac{L}{h};~~~~~\overline{P}=P\times \frac{L}{h}
\end{equation}

In the following, we first present the results of the validation and convergence study and then we present the static response of the fractional-order nonlocal beam for the various loading conditions. For all the numerical results, the mesh discretization was chosen such that $l_f/l_e=10$, where $l_e$ is the length of the discretized element. This choice allows for sufficient number of elements to be included in the horizon of nonlocality at any point hence guaranteeing to accurately capture the nonlocal interactions.

\subsection{Validation}
We validated the performance and evaluated the accuracy of the f-FEM by using the following strategies: \ul{validation \#1}: the f-FEM model was solved for $\alpha=1.0$, corresponding to the classical elastic solid, and the results were compared with those available in the literature\cite{reddy2014introduction}; \ul{validation \#2}: the f-FEM model was solved for non-integer values of $\alpha$ while ignoring the nonlinear effects in Eq.~\eqref{eq: stiffness_mats}, and the results were compared against \cite{patnaik2019FEM}; and \ul{validation \#3}: the f-FEM results were compared with the exact solution for a nonlinear and fractional-order elastic clamped-clamped beam subjected to a distributed transverse load. The validations \#1 and \#2 can rather be considered as checks on the robustness of the nonlinear f-FEM. Specific details of the three validation strategies are provided here below.\\

\begin{figure*}[t!]
    \centering
    \begin{subfigure}[t]{0.49\textwidth}
        \centering
        \includegraphics[width=\textwidth]{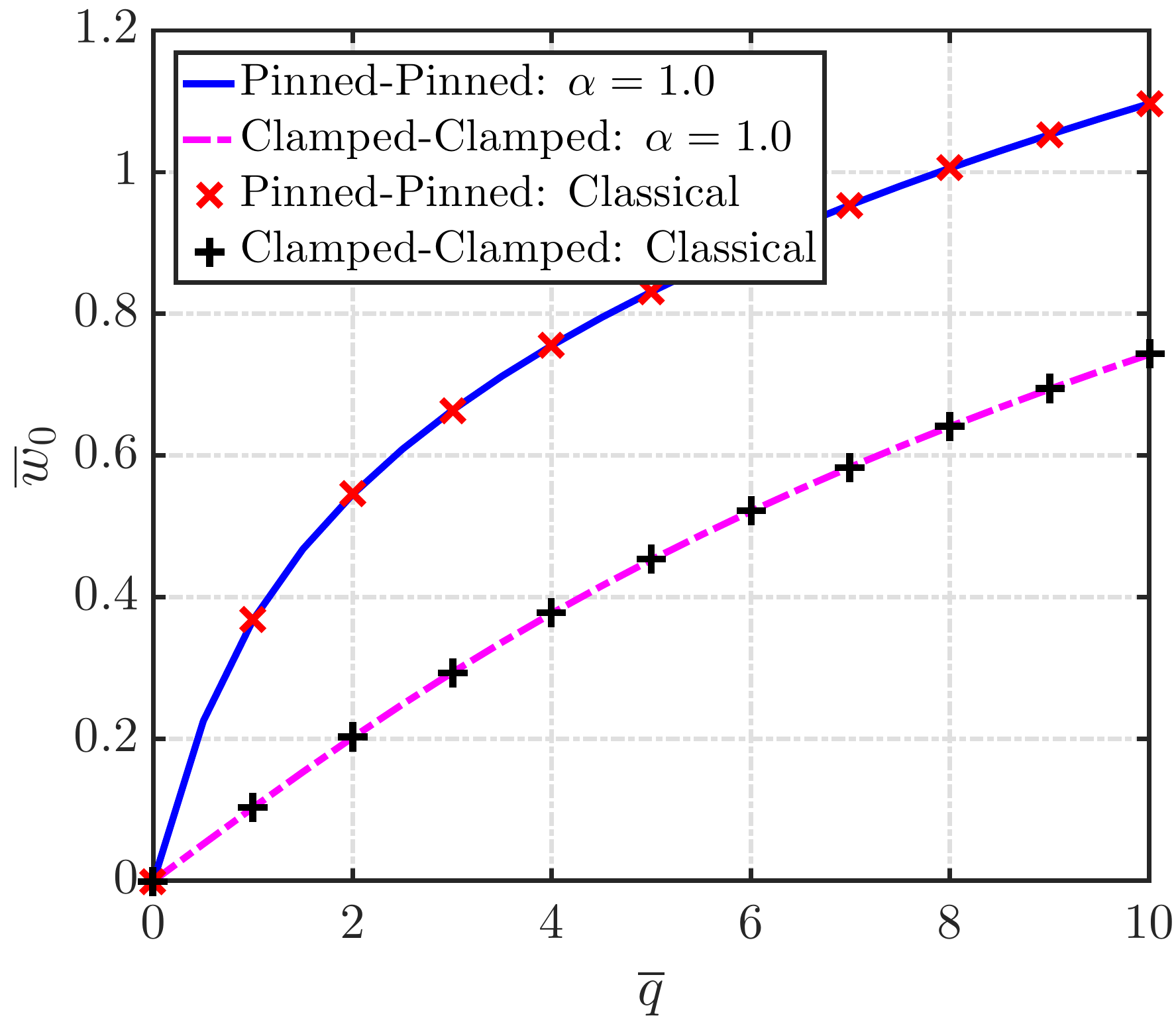}
        \caption{Comparison with the nonlinear local continuum.}
        \label{fig: classical_valid}
    \end{subfigure}%
    ~ 
    \begin{subfigure}[t]{0.49\textwidth}
        \centering
        \includegraphics[width=\textwidth]{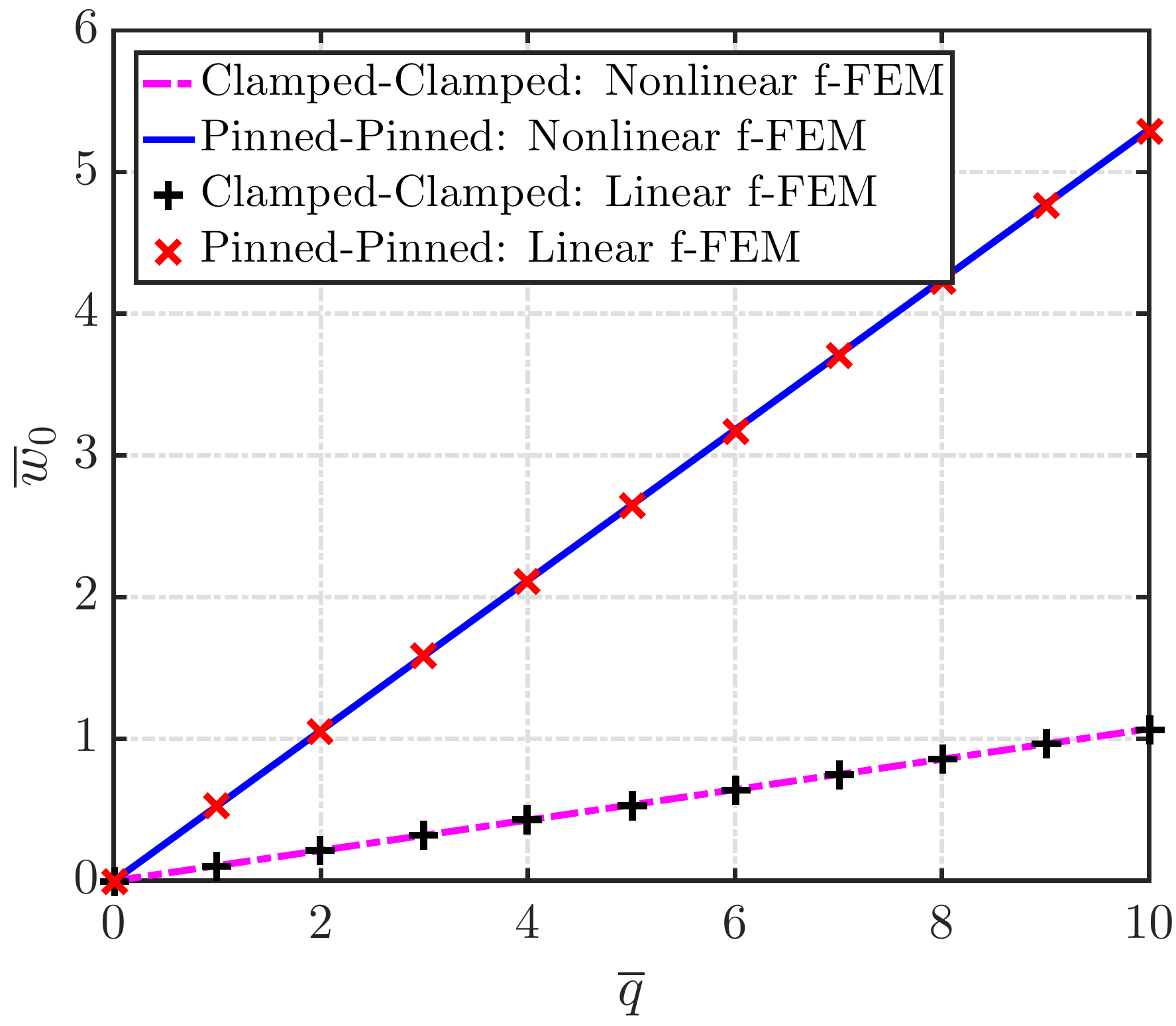}
        \caption{Comparison with the linear nonlocal continnum.}
        \label{fig: linear_valid}
    \end{subfigure}
    \caption{Numerical nonlinear f-FEM results for (a) nonlinear classical continuum $(\alpha=1.0)$ compared with \cite{reddy2014introduction}, and (b) linear nonlocal continuum ($\alpha=0.9$, $l_f=L/10$) compared with \cite{patnaik2019FEM}.}
    \label{fig: valid_2}
\end{figure*}

\noindent
\textbf{Validation \#1}: 
in this study, the fractional order is taken equal to $\alpha=1.0$. This assumption reduces the nonlocal constitutive model (strain-displacement relations) in Eq. \eqref{eq: vonkarman_strain} to their integer-order counterparts for classical elasticity. 
The nonlinear f-FEM is solved for this choice of $\alpha$, $l_f$ being arbitrary for the local elastic solid, and the numerical results are compared with \cite{reddy2014introduction} in the Fig.~(\ref{fig: classical_valid}). This comparison is performed for beams corresponding to two different boundary conditions, namely, pinned-pinned and clamped-clamped. As evident from Fig.~(\ref{fig: classical_valid}), the match between the transverse displacement obtained numerically from the f-FEM and the results in the literature is excellent. The error between them is less than $1\%$ for all cases.\\  

\noindent
\textbf{Validation \#2}: 
in this study, the nonlinear effects are ignored for the sake of a comparison with the literature over fractional-order nonlocal and linear elastic response \cite{patnaik2019FEM}. This simplification is achieved by ignoring the term $D_{x_1}^{\alpha}[w_0(x_1)]$ in the definition for stiffness matrices given in Eq.~\eqref{eq: stiffness_mats}, which contributes to the nonlinearity. The resulting linearized f-FEM model was solved for $\alpha=0.9$ and $l_f=L/10$ and the numerical results are compared with \cite{patnaik2019FEM} in the Fig.~(\ref{fig: linear_valid}). Similar to validation \#1, this comparison is also performed for two different boundary conditions. As evident from Fig.~(\ref{fig: linear_valid}), the match between the transverse displacement obtained numerically from the f-FEM and the results available in the literature is excellent and the relative error is less than $1\%$ for all cases.\\

\begin{figure*}[t!]
    \centering
    \begin{subfigure}[t]{0.49\textwidth}
        \centering
        \includegraphics[width=\textwidth]{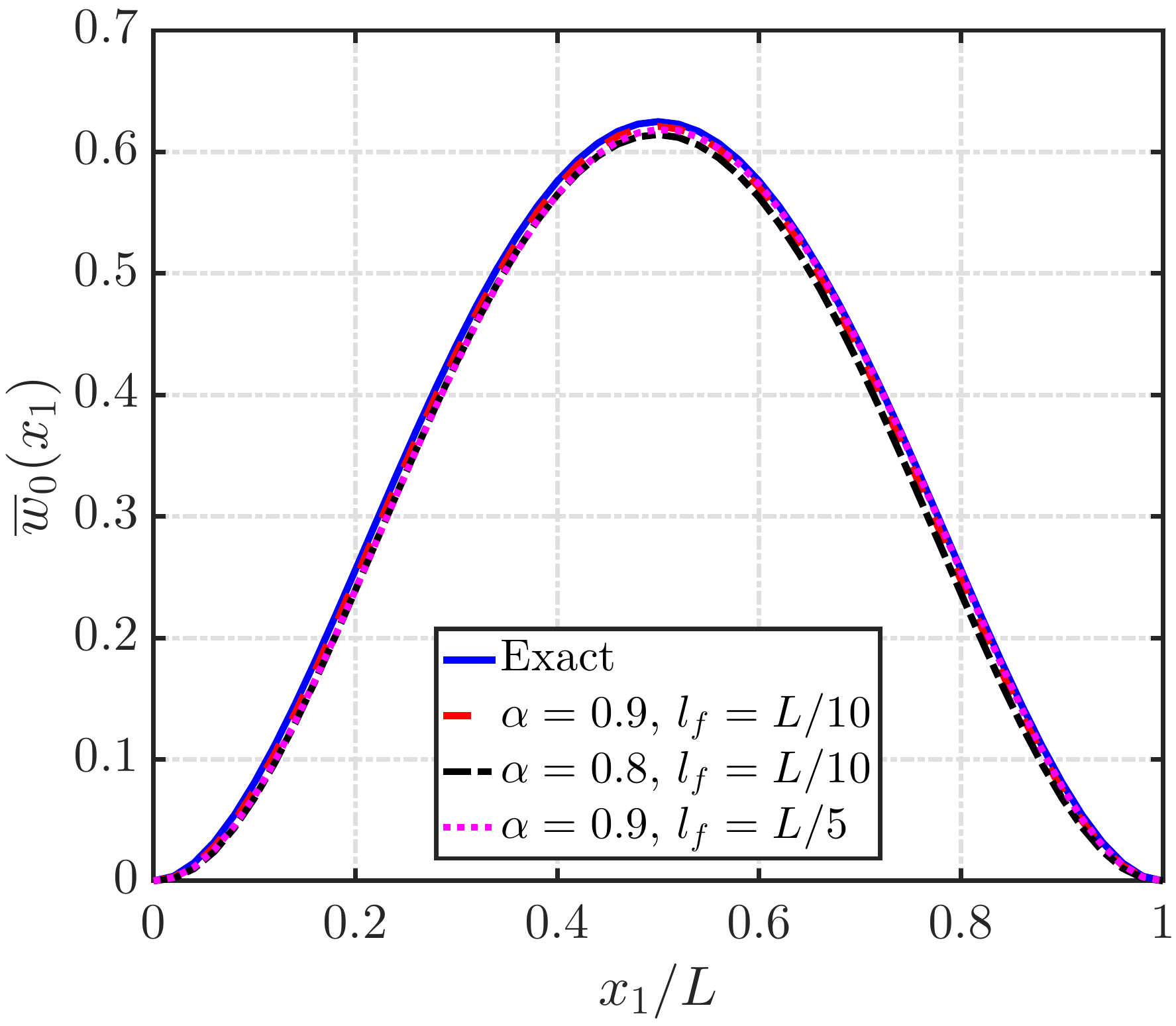}
        \caption{Case 1.}
        \label{fig: validation_cc_1}
    \end{subfigure}%
    ~ 
    \begin{subfigure}[t]{0.49\textwidth}
        \centering
        \includegraphics[width=\textwidth]{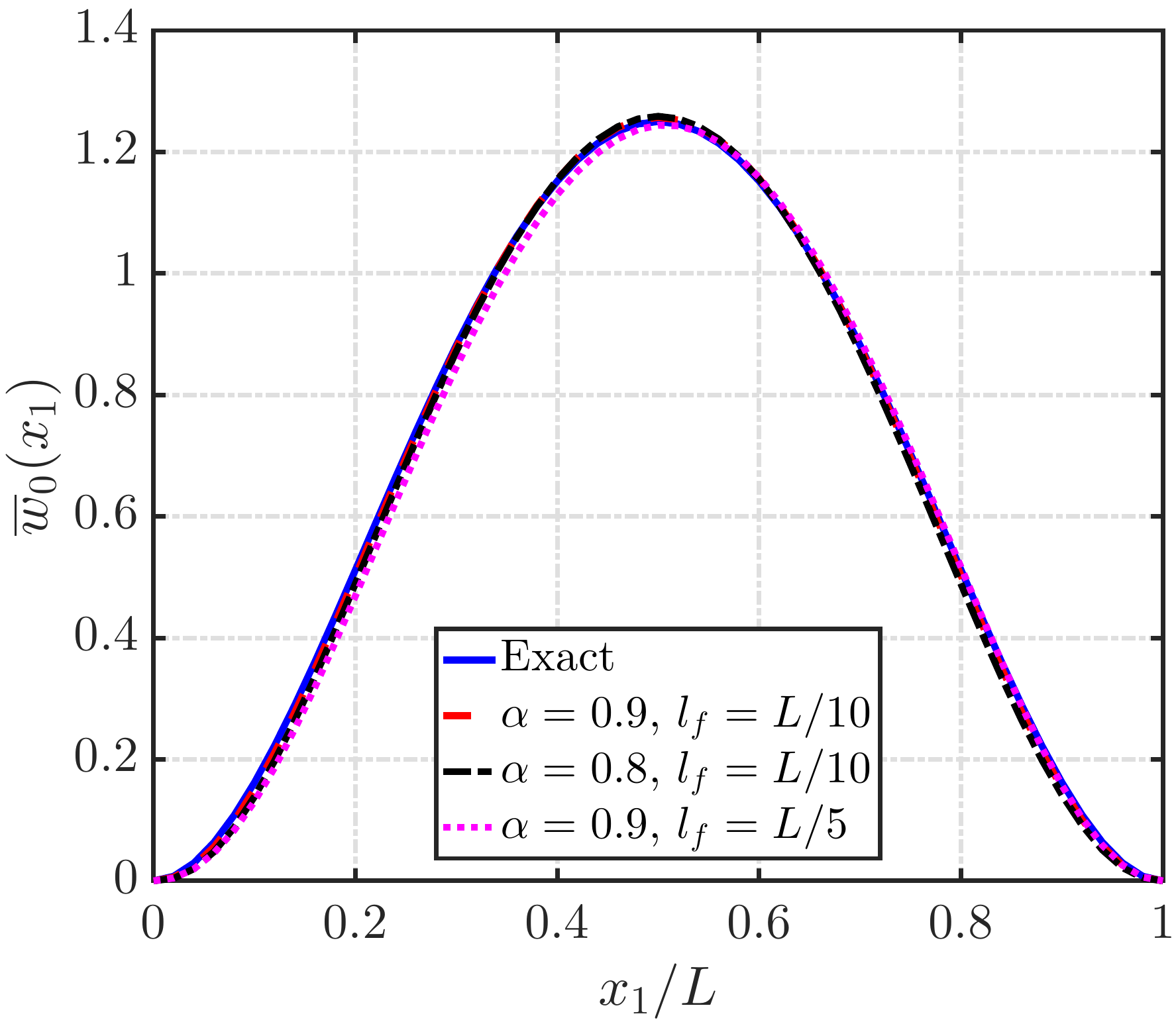}
       \caption{Case 2.}
        \label{fig: validation_cc_2}
    \end{subfigure}
    \caption{Numerical f-FEM results (a) Case 1: $W_0=10h$ and $U_0=0.05h$, and (b) Case 2: $W_0=20h$ and $U_0=0.1h$, corresponding to the assumed displacement field in Eq.~\eqref{eq: valid_cc_disp} compared with the exact solutions for the clamped-clamped boundary condition. Note that, in this validation procedure, the transverse displacement field is assumed independent of $\alpha$ and $l_f$.}
    \label{fig: validation_cc}
\end{figure*}

\noindent
\textbf{Validation \#3}: in this study, the following displacement field at the mid-plane of the beam is assumed:
\begin{subequations}
\label{eq: valid_cc_disp}
\begin{equation}
    u_0(x_1)=U_0\left(1-\frac{x_1}{L}\right)\left(\frac{x_1}{L}\right)
\end{equation}
\begin{equation}
    w_0(x_1)=W_0\left(1-\frac{x_1}{L}\right)^2\left(\frac{x_1}{L}\right)^2
\end{equation}
\end{subequations}
Note that the above displacement fields satisfy the clamped-clamped boundary conditions for the beam given in Eq. \eqref{eq: all_governing_equations}. Although, the above displacement fields are independent of the fractional constitutive parameters $\alpha$ and $l_f$, the axial and the transverse forces required to generate the above elastic response of the beam would be functions of these parameters. These distributed loads are obtained from the strong form of the governing equations given in Eq. \eqref{eq: all_governing_equations} using the constitutive relations stated in Eqs. \eqref{eq: constt_isot} and \eqref{eq: stress_resultants}. Note that these expressions derived from the governing equations contain the appropriate information regarding the geometry. Moreover, the amplitude of the load is dependent on the assumed displacement field distribution in Eq. \eqref{eq: valid_cc_disp}. The closed-form expressions of the distributed axial and transverse loads derived for the exact solutions in Eq. \eqref{eq: valid_cc_disp} are used in the f-FEM model given in Eq. \eqref{eq: algebraic_eqs}. This strategy is similar to the one presented in \cite{patnaik2019FEM}. The numerical results obtained from the f-FEM for the mid-plane transverse displacement $w_0(x_1)$ are compared with the exact solution. The results, shown in Fig.~(\ref{fig: validation_cc}), are in excellent agreement. The error between the numerically obtained f-FEM solutions and the exact solutions is less than $3\%$ for all cases. Similar inverse approaches for validation have also been adopted in \cite{zhuang2009numerical,zhang2013novel} for different numerical methods.

\subsection{Convergence}
\label{sec: convergence_results}
In the following, we present the results of the sensitivity analysis on the FE mesh size. The convergence of the integer-order FEM with element discretization, referred to as the $h-$refinement, is well established in the literature \cite{reddy2014introduction}. In this study, we noted that in addition to the FE mesh size, the convergence of the f-FEM also depends on the fractional order $\alpha$, on the relative size of $l_e$, and on the nonlocal length scale $l_f$. It appears that the convergence of the f-FEM depends on the strength of the nonlocal interactions across the nonlocal horizon. Therefore, sufficient number of elements $N_{inf}(=l_f/l_e)$ should be available in the influence zone of the fractional-order (nonlocal) interactions. Additionally, mesh refinement across the length of the 1D beam would increase the spatial resolution and hence, result in lower inconsistencies due to the truncation of the nonlocal horizon caused by the ceil and floor operations connected to Eq.~(\ref{eq: b_mat_int_Scheme}). This approach would provide a better numerical approximation for the nonlocal matrices, $[\tilde{B}_\square]$ in Eqs. \eqref{eq: b_mat_numer_simp1}-\eqref{eq: b_mat_int_Scheme}. Therefore, convergence of the f-FEM is expected when the number of elements in the influence zone $N_{inf}$, referred to as the “dynamic rate of convergence” \cite{norouzzadeh2017finite,patnaik2019FEM}, is sufficient to accurately capture the nonlocal interactions. 

\begin{figure*}[h!]
    \centering
    \begin{subfigure}[t]{0.49\textwidth}
        \centering
        \includegraphics[width=\textwidth]{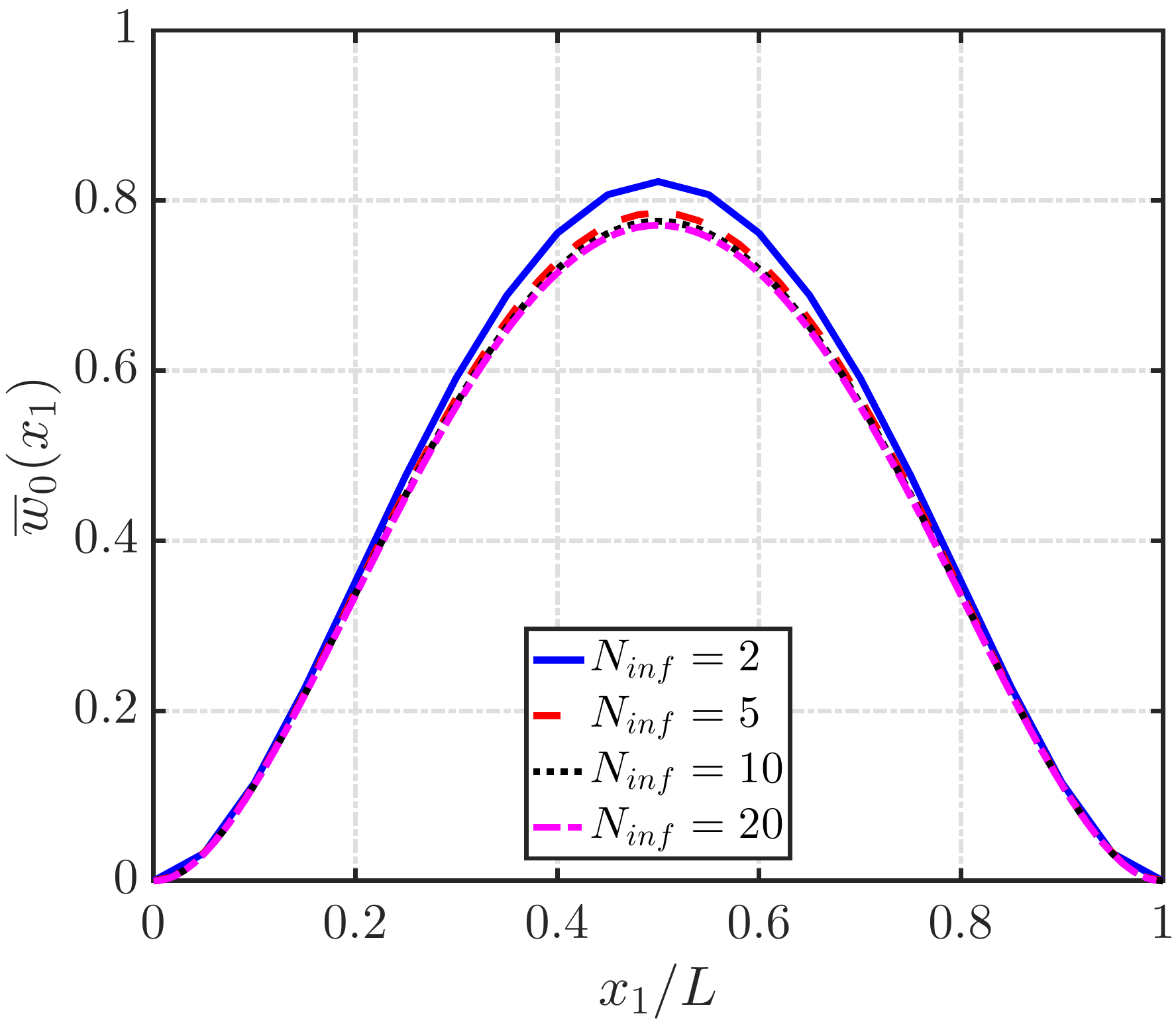}
        \caption{Clamped-Clamped beam.}
        \label{fig: convergence_cc}
    \end{subfigure}%
    ~ 
    \begin{subfigure}[t]{0.49\textwidth}
        \centering
        \includegraphics[width=\textwidth]{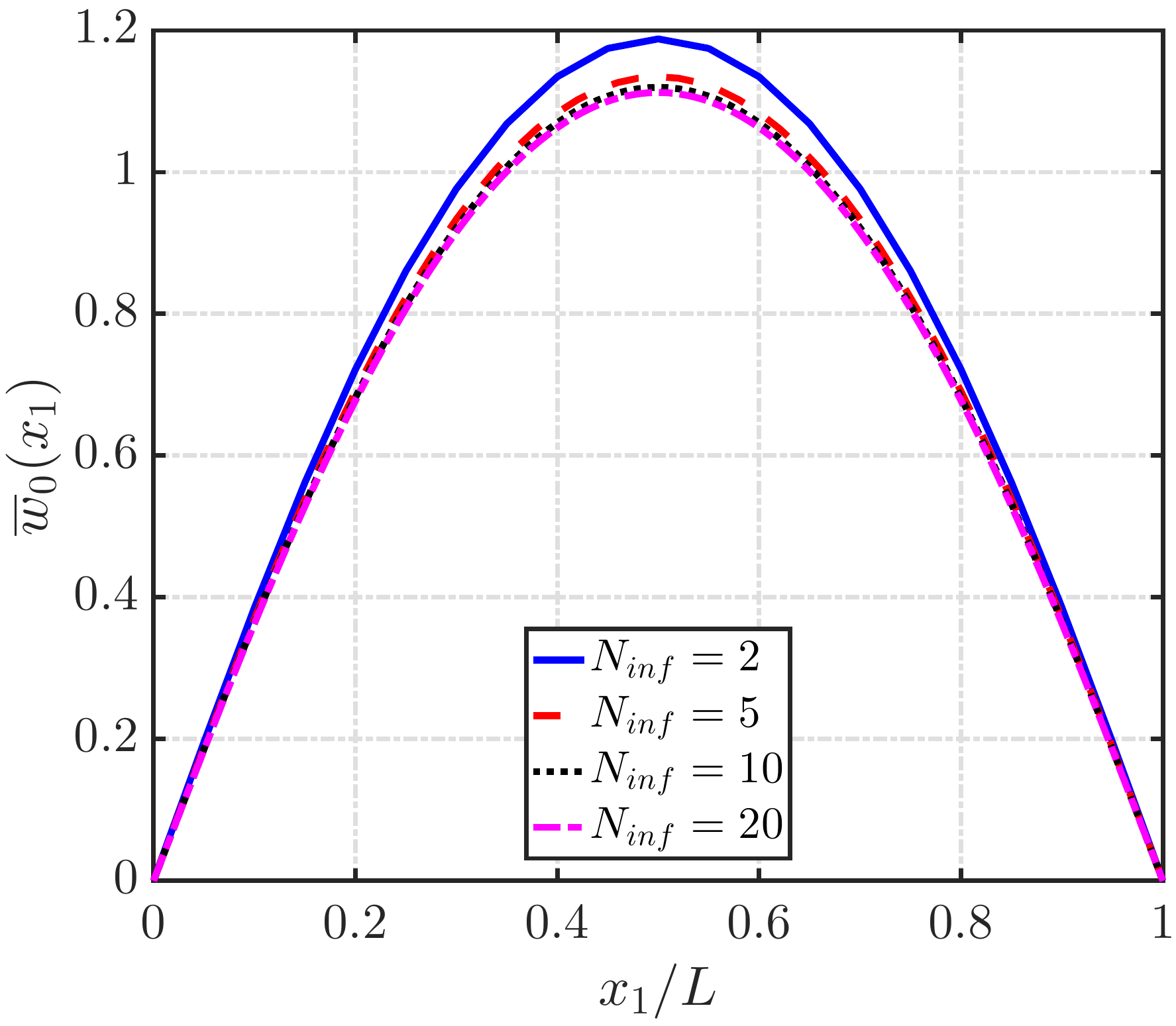}
       \caption{Pinned-Pinned beam.}
        \label{fig: convergence_ss}
    \end{subfigure}
    \caption{Numerical convergence of the f-FEM numerical model ($\alpha=0.8$ and $l_f=L/10$).}
    \label{fig: convergence_cases}
\end{figure*}

In the context of the current study, we establish the convergence of the f-FEM for two different boundary conditions. The fractional-order parameters are chosen as $\alpha=0.8$ and $l_f=L/10$. The results are presented in Fig.~(\ref{fig: convergence_cases}). In both cases, Fig.~(\ref{fig: convergence_cc}) for the clamped-clamped and Fig.~(\ref{fig: convergence_ss}) for the pinned-pinned beams, the dynamic rate of convergence is noted to be $N_{inf}\geq 5$. Additionally, a comprehensive comparison of the convergence of the f-FEM for various values of the fractional parameters $\alpha$ and $l_f$, for a clamped-clamped beam, is presented in Table \ref{tab: convergence_comp}. We have discussed in \S\ref{ssec: effect}, the reason behind choosing a restricted range of $\alpha \in [0.5,1]$. In Table \ref{tab: convergence_comp}, the convergence may be assessed by comparing the values along a prescribed column on the far right, corresponding to the transverse displacement for a chosen set of fractional constitutive parameters. As evident from Table \ref{tab: convergence_comp}, for higher values of $l_f$ and lower values of $\alpha$, that is for stronger fractional-order nonlocal interactions, larger number of elements are necessary for the convergence of f-FEM. Note also that $N_{inf}=10$ provides a satisfactory numerical convergence yielding a difference smaller than $1\%$ between successive refinements. Therefore, the f-FEM simulations presented below will be carried out using this mesh discretization. 

\begin{table}[h!]
    \centering
    \begin{tabular}{c | c |c c c c c c}
    \hline\hline
       \multirow{2}{14em}{ Fractional Horizon Length, $l_f$} & \multirow{2}{6em}{~~~~~~$N_{inf}$} & \multicolumn{6}{c}{fractional-order, $\alpha$}\\
         \cline{3-8}
         &  & $\alpha=1.0$& $\alpha=0.9$& $\alpha=0.8$ & $\alpha=0.7$& $\alpha=0.6$& $\alpha=0.5$  \\
         \hline\hline
          \multirow{4}{5em}{$l_f=L/5$} & 2 & 0.7352 & 0.7911 & 0.8454 & 0.9004 & 0.9612 & 1.0367\\
          & 5 & 0.7394 & 0.7819 & 0.8202 & 0.8569 & 0.8967 & 0.9488\\
          & 10 & 0.7426 & 0.7821& 0.8168 & 0.8494 & 0.8844 & 0.9322\\
          & 20 & 0.7428 & 0.7810 & 0.8140 & 0.8449 & 0.8775 & 0.9224\\
          \hline
         \multirow{4}{5em}{$l_f=L/10$} & 2 & 0.7410 & 0.7821 & 0.8223 & 0.8636 & 0.9095 & 0.9667\\
         & 5 & 0.7426 & 0.7653 & 0.7856 & 0.8054 & 0.8273 & 0.8556\\
         & 10 & 0.7428 & 0.7606 & 0.7756 & 0.7899 & 0.8057 & 0.8264\\
         & 20 & 0.7429 & 0.7583 & 0.7708 & 0.7826& 0.7956 & 0.8128\\
         \hline
          \multirow{4}{5em}{$l_f=L/20$} & 2 & 0.7424 & 0.7807  & 0.8187 & 0.8577 & 0.8994 & 0.9472\\
          & 5 & 0.7428 & 0.7597 & 0.7748 & 0.7895 & 0.8053 & 0.8244\\
          & 10 & 0.7426 & 0.7538 & 0.7627 & 0.7710 & 0.7802 & 0.7920\\
          & 20 & 0.7429 & 0.7510 & 0.7568 & 0.7622 & 0.7684 & 0.7769\\
          \hline\hline
    \end{tabular}
    \caption{Convergence of the f-FEM for clamped-clamped beam with different fractional parameters.}
    \label{tab: convergence_comp}
\end{table}


\subsection{Effect of the fractional-order constitutive parameters}
\label{ssec: effect}
In this section, we use the f-FEM to analyze the static response of the fractional-order nonlocal and geometrically nonlinear Euler-Bernoulli beam. More specifically, we analyze the effect of the fractional parameters $\alpha$ and $l_f$ on the response of the beam subject to different loading and boundary conditions.

\begin{figure*}[ht!]
    \centering
    \begin{subfigure}[t]{0.49\textwidth}
        \centering
        \includegraphics[width=\textwidth]{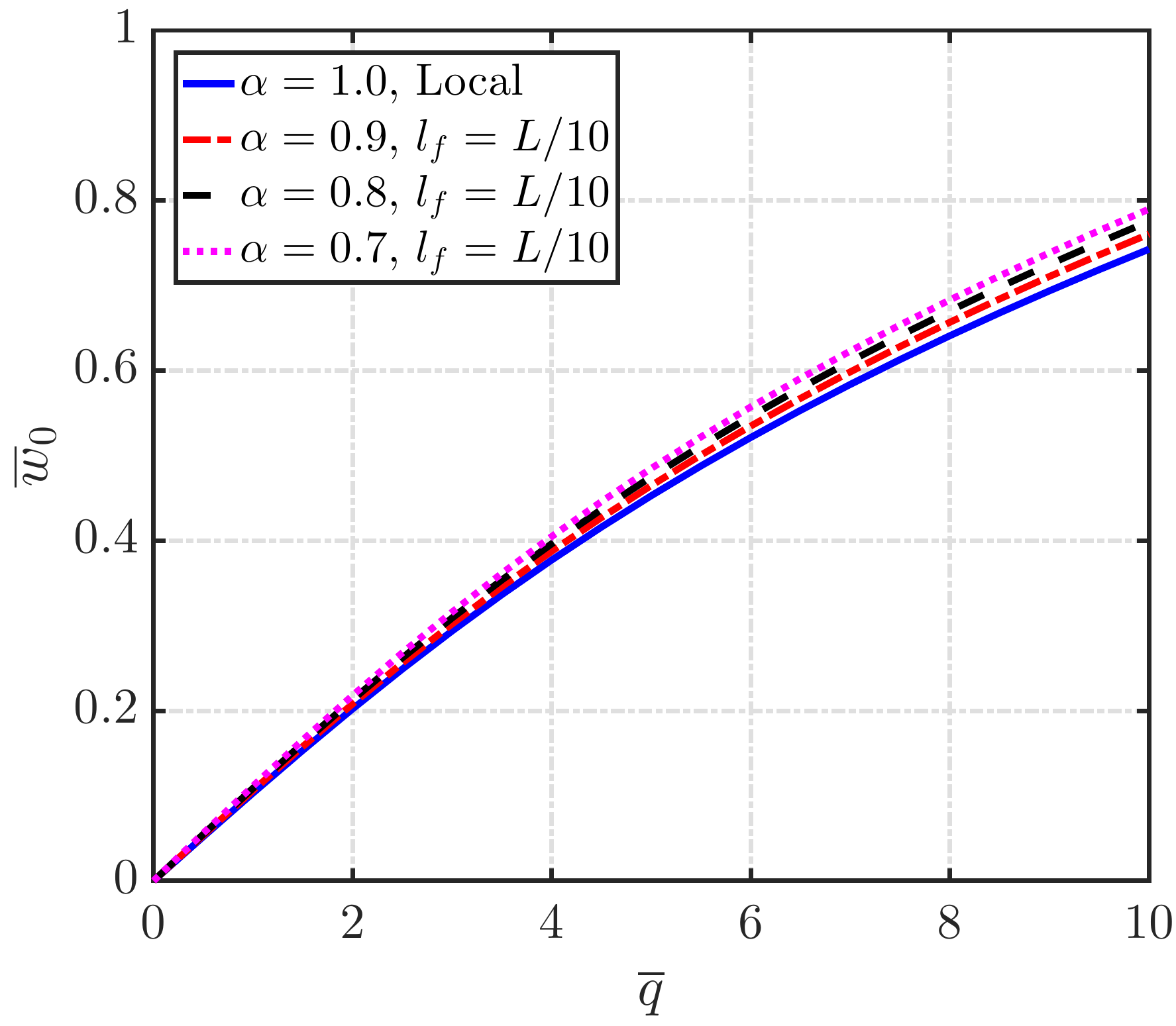}
        \caption{$\overline{w}$ vs $\alpha$.}
        \label{fig: cc_w_alpha}
    \end{subfigure}%
    ~ 
    \begin{subfigure}[t]{0.49\textwidth}
        \centering
        \includegraphics[width=\textwidth]{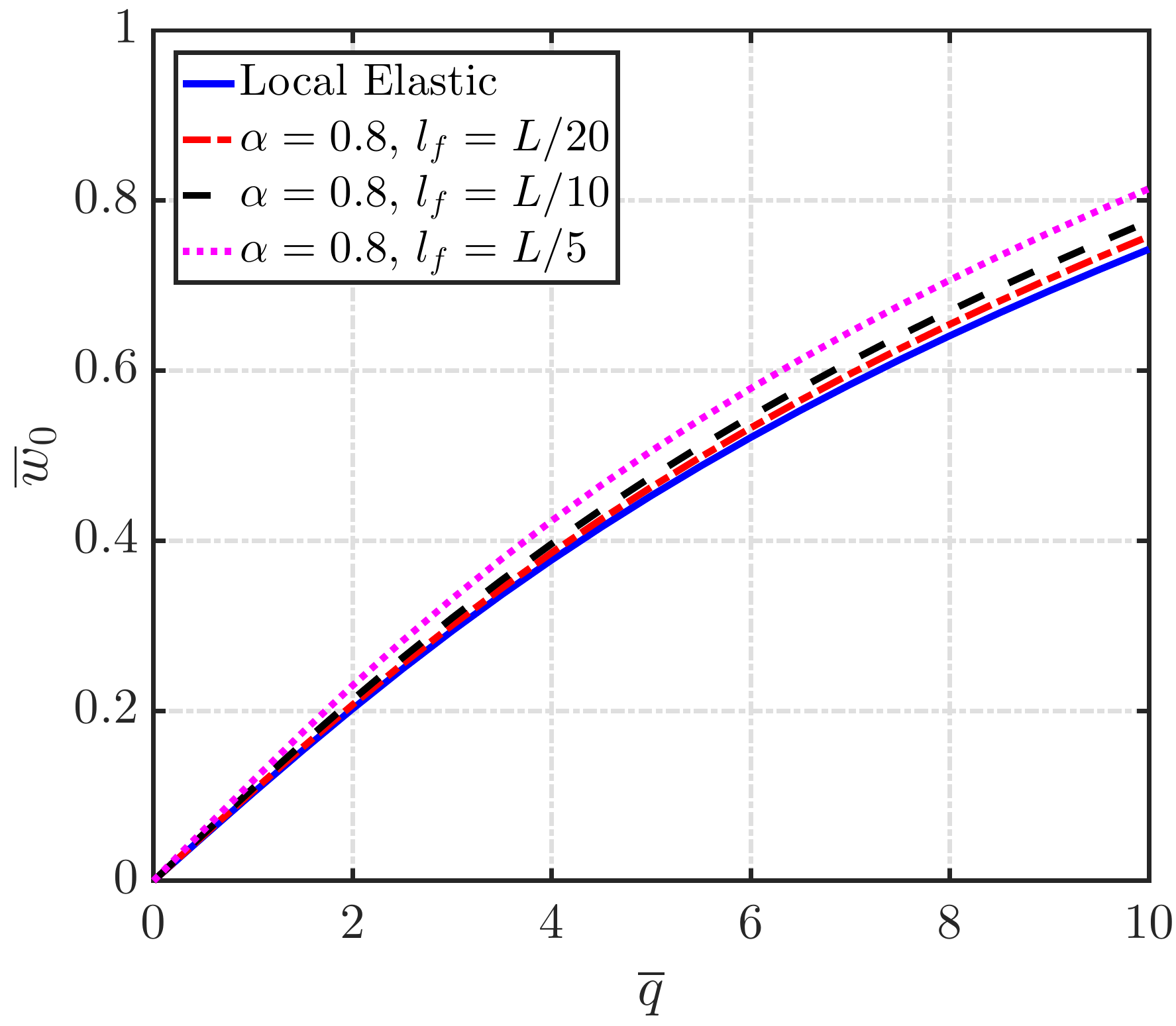}
        \caption{$\overline{w}$ vs $l_f$.}
        \label{fig: cc_w_lf}
    \end{subfigure}
    \caption{Nonlocal elastic transverse displacement of the clamped-clamped beam compared for different values of the fractional-order constitutive properties.}
    \label{fig: cc_w}
\end{figure*}

We start by considering a beam that is clamped at both ends and is subject to a UDL of magnitude $q_0$ (in N/m). The nonlinear load-displacement curves of the clamped-clamped beam are obtained and compared for different values of $\alpha$ and $l_f$. The results of this comparison are provided in Fig.~(\ref{fig: cc_w}). The effect of the fractional order $\alpha$, with $l_f$ being held constant, is studied in Fig.~(\ref{fig: cc_w_alpha}). Recall that for $\alpha=1.0$, the differential and algebraic form of the governing equations and boundary conditions, reduce to those of a geometrically nonlinear classical (local) elasticity study.
Additionally, we also compare the effect of the horizon of nonlocality over the elastic response of the beam in Fig.~(\ref{fig: cc_w_lf}). As evident from Fig.~(\ref{fig: cc_w}), the transverse displacement increases with increasing degree of nonlocality, that is by reducing the value of $\alpha$ or by increasing the horizon length $l_f/L$. These trends point to a reduction in the stiffness of the fractional-order nonlocal structure. As expected, a convergence of the fractional-order model to the local elastic response is noted for $\alpha$ close to $1.0$ and $l_f/L<<1$.

\begin{figure*}[ht!]
    \centering
    \begin{subfigure}[t]{0.49\textwidth}
        \centering
        \includegraphics[width=\textwidth]{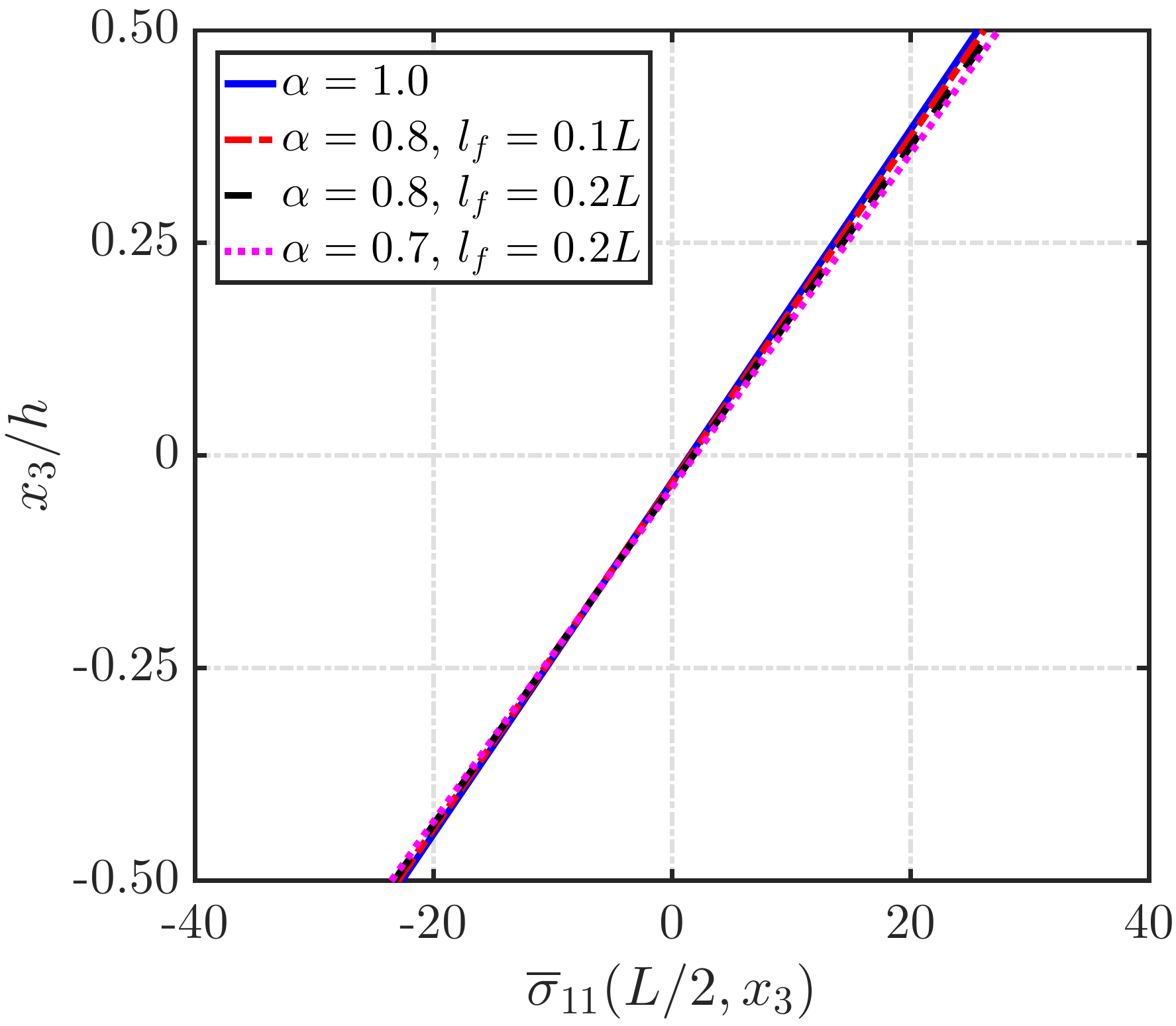}
        \caption{$\overline{q}=2$.}
        \label{fig: cc_lin_signa}
    \end{subfigure}%
    ~ 
    \begin{subfigure}[t]{0.49\textwidth}
        \centering
        \includegraphics[width=\textwidth]{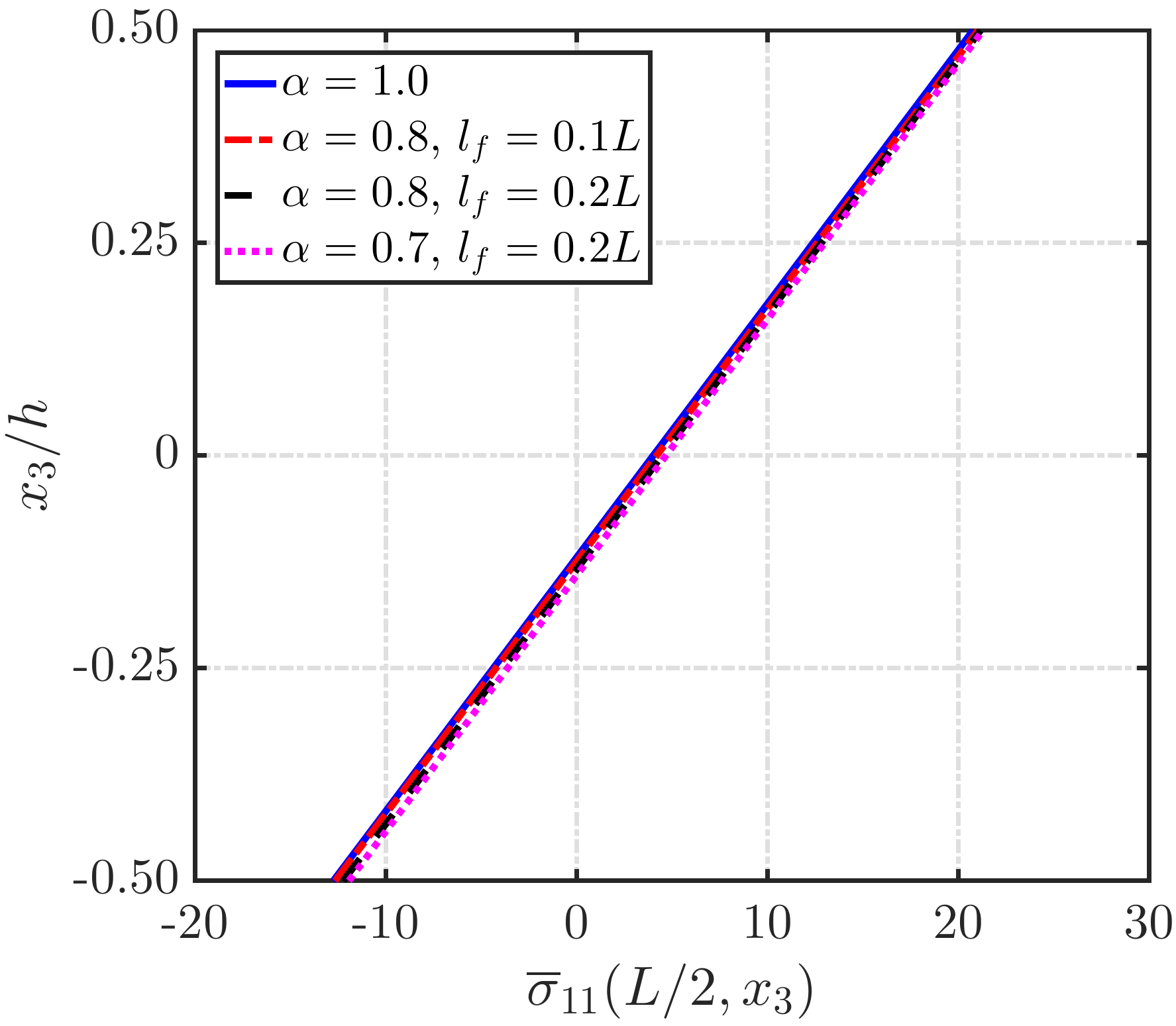}
        \caption{$\overline{q}=10$.}
        \label{fig: cc_nl_sigma}
    \end{subfigure}
    \caption{Normal axial stress $\tilde{\sigma}_{11}$ in the nonlocal beam subject to clamped-clamped boundary conditions and two cases of the transverse load compared for different fractional-order constitutive properties.}
    \label{fig: cc_sig}
\end{figure*}

For the clamped-clamped beam described above, we also analyzed the effect of the fractional constitutive parameters on the axial normal stress across the thickness of the beam at its mid-length, that is at $x_1=L/2$. The results are presented in Fig.~(\ref{fig: cc_sig}) in terms of the non-dimensional stress following the application of two different transverse loads. As expected of nonlinear response, an increase in magnitude of the transverse load (for $\bar{q}=10)$ produces a non-zero axial stress at the mid-plane (see Fig.~(\ref{fig: cc_nl_sigma})) of the beam. This axial stress caused by the transverse loading is a direct result of the axial-transverse coupling induced by the geometric nonlinearity. For lower magnitudes of the transverse load $\overline{q}=2$, a weak nonlinearity is clear from the lower value of $\bar{\sigma}_{11}$ at the mid-plane (see Fig.~(\ref{fig: cc_lin_signa})). Note that, given the linear nature of the constitutive relations in Eq. \eqref{eq: constt_isot} used for the current analysis, the observations noted above for the stress in the beam, directly extend to the fractional-order strain.

Before proceeding further, we briefly discuss how the degree of nonlocality increases by reducing the order $\alpha$ and increasing the length-scale $l_f$. Note that the convolution kernel used within the definition of the fractional-order derivative is a power-law kernel which is monotonically decreasing in nature. It appears that, by reducing the value of the order $\alpha$, the magnitude of this kernel increases. More specifically, for a fixed point $x$ interacting with a specific point $s$ in its horizon of nonlocality, $\left(\frac{1}{|x-s|}\right)^{\alpha_1} > \left(\frac{1}{|x-s|}\right)^{\alpha_2} \forall~\alpha_1 < \alpha_2$. Since the magnitude of the kernel increases, the strength of the corresponding nonlocal interactions increases and so does the degree of nonlocality. Similarly, by increasing the length scale $l_f$, the size of the horizon of nonlocality increases. Hence, information corresponding to the nonlocal interactions between a larger number of points within the solid are accounted for in the formulation. Thus the degree of nonlocality increases with increasing values of $l_f$.

\begin{figure*}[ht!]
    \centering
    \begin{subfigure}[t]{0.49\textwidth}
        \centering
        \includegraphics[width=\textwidth]{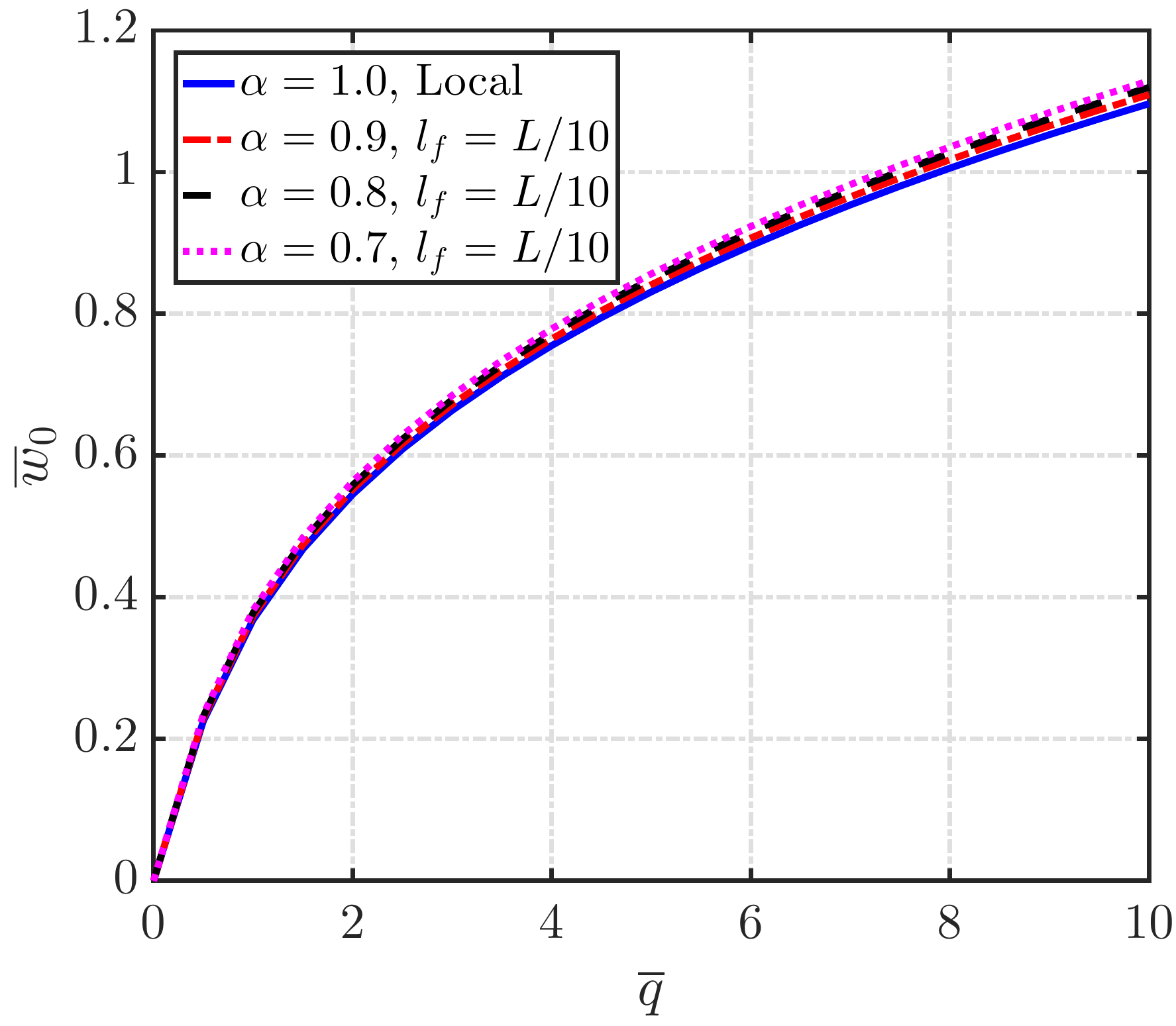}
        \caption{$\overline{w}$ vs $\alpha$.}
        \label{fig: ss_w_alpha}
    \end{subfigure}%
    ~ 
    \begin{subfigure}[t]{0.49\textwidth}
        \centering
        \includegraphics[width=\textwidth]{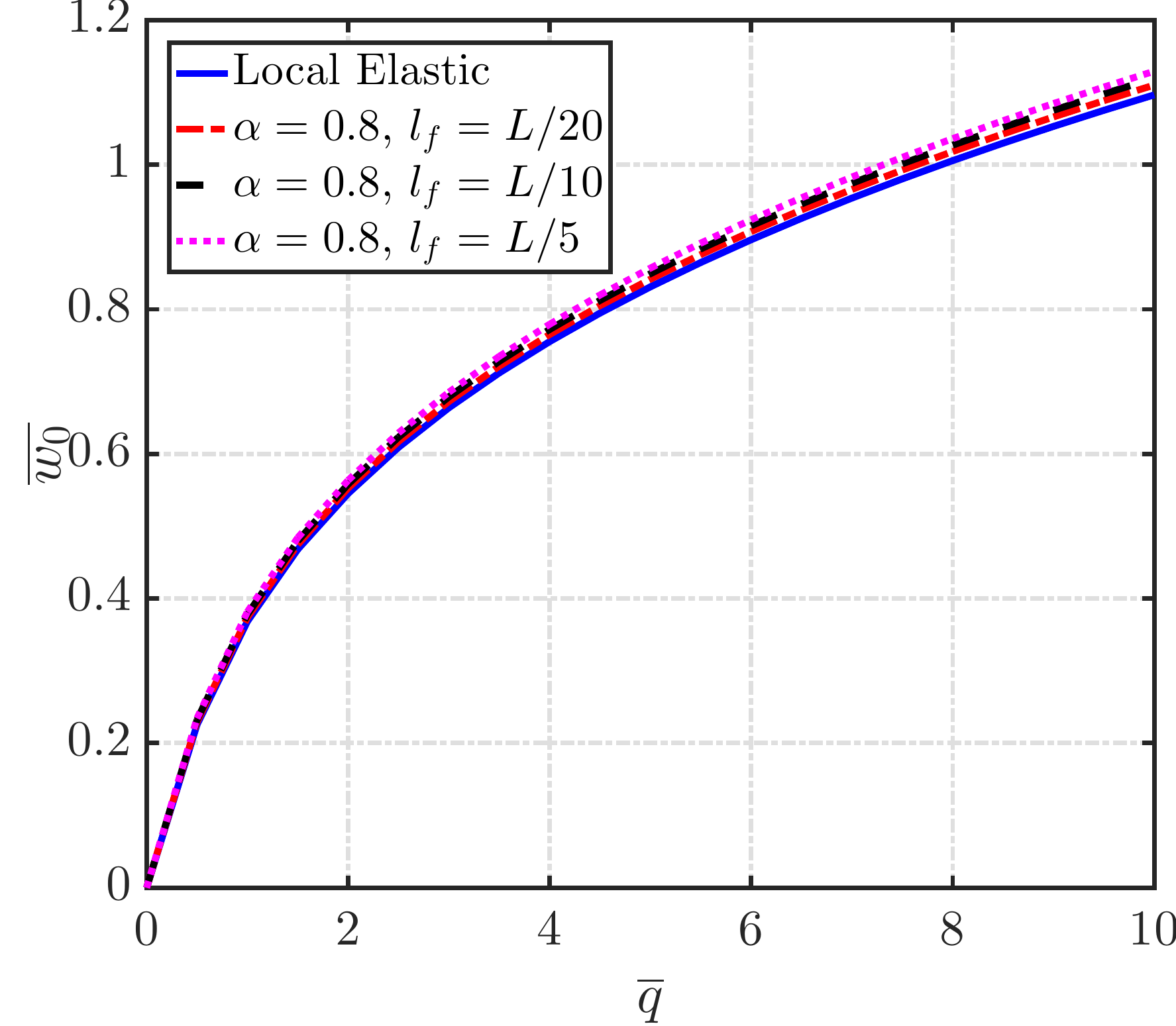}
        \caption{$\overline{w}$ vs $l_f$.}
        \label{fig: ss_w_lf}
    \end{subfigure}
    \caption{Comparison of the nonlocal elastic transverse displacement of the pinned-pinned beam for different values of the fractional-order constitutive properties.}
    \label{fig: ss_w}
\end{figure*}

We extend the above study to analyze a beam with both ends pinned. The effect of the fractional model parameters over the nonlocal elastic response of the geometrically nonlinear pinned-pinned beam are shown in Fig.~(\ref{fig: ss_w}). From these figures, it may be noted that the transverse displacement increases while reducing the value of $\alpha$ and increasing the nonlocal horizon length $l_f$. These observations for the pinned-pinned beam agree very well with those noted previously for the clamped-clamped beam. The consistency of these observations regardless of the nature of the boundary conditions indicates the robustness of the fractional order model and marks a net departure from the paradoxical results obtained in the literature for more traditional integral beam theories \cite{challamel2008small,khodabakhshi2015unified,fernandez2016bending}. These inaccurate predictions have been traced back to the ill-posed nature of the nonlocal constitutive relations as discussed in detail in \cite{romano2017constitutive}. On the contrary, the fractional-order model of nonlocality is free from mathematical ill-posedness and leads to consistent predictions of softening for the nonlocal beam structures \cite{patnaik2019FEM}. 

\begin{figure*}[ht!] 
    \centering
    \begin{subfigure}[t]{0.49\textwidth}
        \centering
        \includegraphics[width=\textwidth]{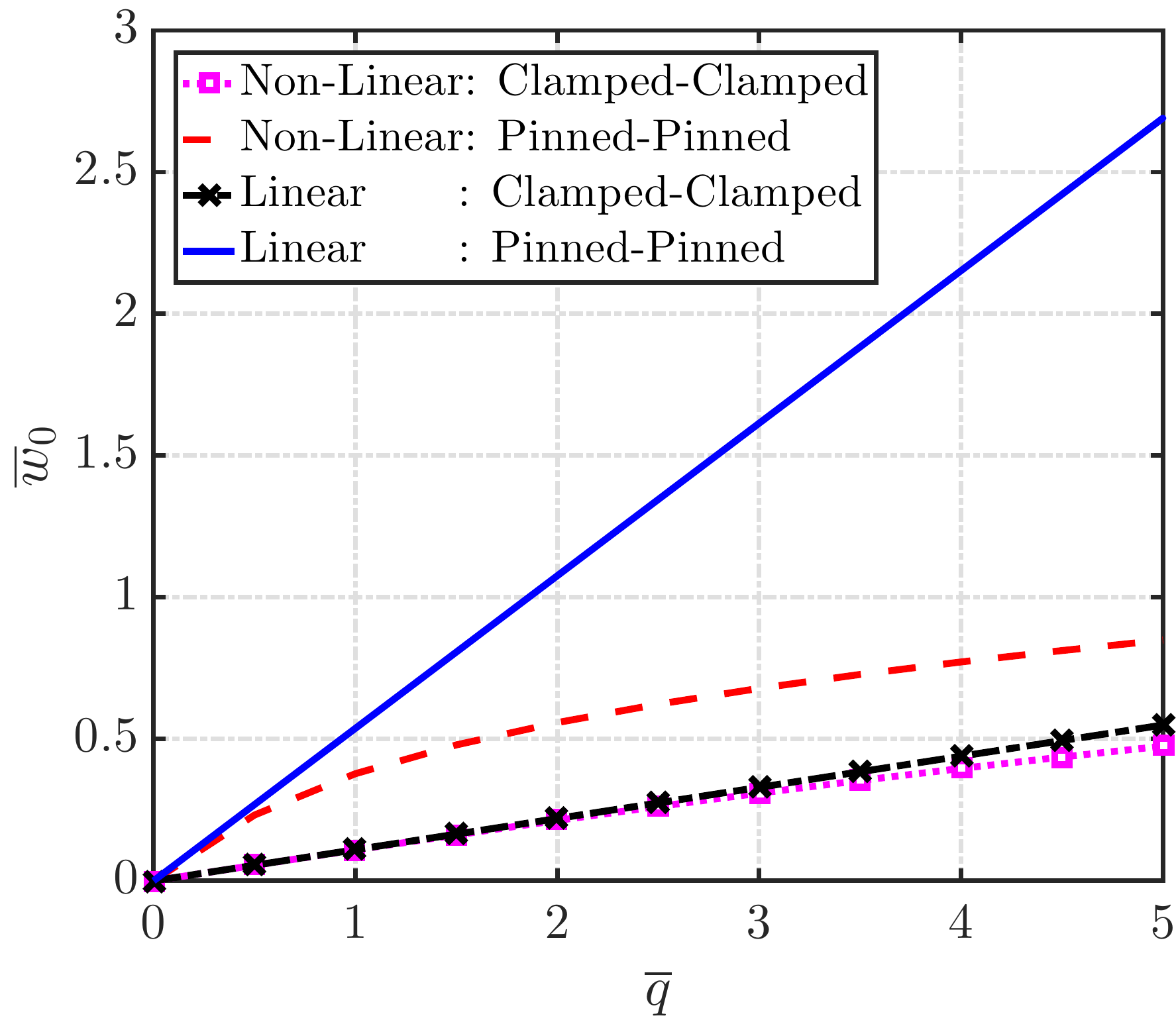}
        \caption{UDL, $q_0$ per unit length.}
        \label{fig: cc_vs_ss_udl}
    \end{subfigure}%
    ~ 
    \begin{subfigure}[t]{0.49\textwidth}
        \centering
        \includegraphics[width=\textwidth]{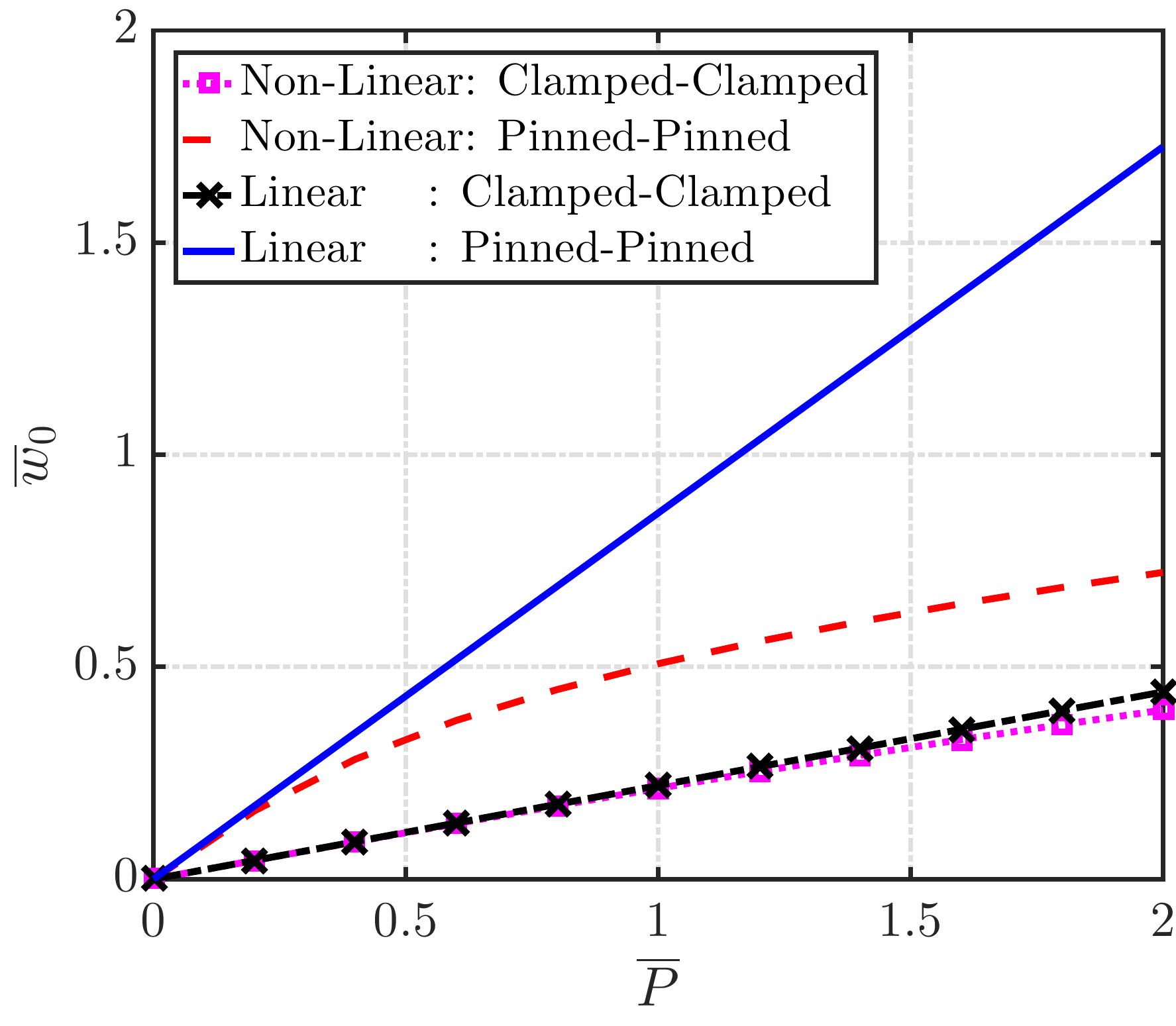}
        \caption{PtL, $P$.}
        \label{fig: cc_vs_ss_ptl}
    \end{subfigure}
    \caption{Geometric non-linear effects compared for loading and boundary conditions ($\alpha=0.8$, $l_f=L/10$).}
    \label{fig: cc_vs_ss}
\end{figure*}

Finally, the nonlocal response of the clamped-clamped and pinned-pinned beams are compared for different loading conditions and the results are shown in Fig.~(\ref{fig: cc_vs_ss}). The linear and nonlinear load-displacement curves of the Euler-Bernoulli beam subject to an UDL and with fractional constitutive parameters fixed at $\alpha=0.8$ and $l_f=L/10$ are compared in Fig.~(\ref{fig: cc_vs_ss_udl}). A similar comparison for the beam subject to a PtL applied at mid-length is shown in Fig.~(\ref{fig: cc_vs_ss_ptl}). It appears that the effect of the geometric nonlinearity is significant for the nonlocal pinned-pinned beams, irrespective of the loading conditions. Note that this is in agreement with the observations from the classical theory of elasticity \cite{reddy2014introduction}.  

An important remark should be made concerning the physically acceptable range for the order $\alpha$. Although a consistent softening is observed for decreasing value of $\alpha$, we highlight here that for very low values of $\alpha$, non-physical results were noted in literature\cite{sumelka2014thermoelasticity,sumelka2015non}. Note that the fractional-order kinematic relations of order $\alpha\in(0,1)$ define the strain as a fractional-order interpolation of the zeroth and first-order derivatives of the displacement. Hence, there exists a limit on the order of the fractional derivative beyond which this fractional model of nonlocality breaks down. In the nonlinear study, we note that this specific value is $\alpha \approx 0.4$. It was also observed that this critical value for $\alpha$ depends on other parameters including the relative ratio of $l_f/L$, where $l_f$ is the nonlocal length scale and $L$ is the geometric length of the structure. However, we emphasize that this breakdown is not a characteristic of the f-FEM technique employed for the numerical analysis, but instead a physical limit of the fractional-order modeling\cite{patnaik2019FEM}.


\section{Conclusions}
\label{sec: conclusion}
The current study laid the foundation for the analysis of geometrically nonlinear and nonlocal solids via fractional-order operators while providing accurate numerical schemes for their solution. More specifically, the geometrically nonlinear formulation for fractional-order nonlocal Euler-Bernoulli beams was developed and its solution was achieved via a dedicated finite-element-based numerical method. The general fractional-order continuum theory was extended to develop a frame-invariant and dimensionally consistent nonlinear fractional-order strain-displacement theoretical framework. This theory is well suited to capture scale effects, nonlocality, and heterogeneity in both complex solids and interfaces undergoing large deformations. The constitutive model for the nonlinear response of fractional-order beams was derived following variational principles. The 1D beam model lays the foundation for a more general nonlinear analysis of fractional-order nonlocal solids. Owing to the complexity of the nonlocal models, a fractional-order FEM (f-FEM) was also developed by means of potential energy minimization. This nonlinear f-FEM provides an efficient and accurate computational approach to obtain the numerical solution of nonlinear integro-differential governing equations with a singular kernel. The robustness of this numerical approach for the solution of nonlocal elasticity problems was demonstrated by showing the consistency of the predicted nonlocal response subject to different loads and boundary conditions. This latter aspect is in contrast with other traditional approaches to nonlocal elasticity that, occasionally, provide paradoxical results depending on the selection of both loading and boundary conditions. Finally we note that, although the f-FEM was used in this study to solve nonlinear beam models, the same computational framework can be easily extended either to the simulation of systems with higher dimensions or to other nonlinear physical processes described by fractional differ-integral equations.

\section*{Appendix 1}
\textbf{Frame invariance of the fractional-order kinematic relations: }
To establish frame invariance of the formulation, we begin with the definition given in Eq.~\eqref{eq: Fractional_F_net} for the deformation gradient tensor following fractional-order kinematics. The frame-invariance of each of the individual terms in the definition would establish the frame invariance of the complete tensor $\tilde{\textbf{F}}_{X}$. For this purpose, we consider a rigid-body motion superimposed on a general point $\textbf{X}$ (see Fig.~(1)) of the reference configuration of the body as:
\begin{equation}
\label{eq: rigid}
\bm{\Psi}(\textbf{X})=\textbf{c}+\textbf{Q}\textbf{X},
\end{equation}
where $\textbf{Q}$, a proper orthogonal tensor, and $\textbf{c}$, a spatially constant term, represent the rotation and translation of the reference frame, respectively. Under this rigid-body motion, the fractional deformation gradient tensor should transform as $\tilde{\textbf{F}}^{\Psi}_X = \textbf{Q}$ (similar to the classical continuum elasticity) such that $\tilde{\textbf{F}}^{\Psi T}_X\tilde{\textbf{F}}^{\Psi}_X = \textbf{I}$. This would result in trivial strain measures for the rigid body transformation. The above rigid body motion employed in the definition of $\tilde{\textbf{F}}_X$ given in Eq.~(\ref{eq: Fractional_Fa}) results in:
\begin{equation}
\label{eq: def_grad}
\begin{split}
\tilde{\textbf{F}}_{X_{ij}}^{\Psi} = \frac{1}{2} \Gamma(2-\alpha) \biggl[ \frac{L_{A_{j}}^{\alpha-1}}{\Gamma(1-\alpha)} \int_{X_{A_{j}}}^{X_j} \frac{D^1_{R_j} \Psi_{i}(\textbf{R})}{(X_j-R_j)^{\alpha}}  \mathrm{d}R_j + \frac{L_{B_{j}}^{\alpha-1}}{\Gamma(1-\alpha)} \int_{X_{j}}^{X_{B_j}} \frac{D^1_{R_j} \Psi_{i}(\textbf{R})}{(R_j-X_j)^{\alpha}} \mathrm{d}R_j \biggr]
\end{split}
\end{equation}
where $\textbf{R}$ is a dummy vector representing the spatial variable of a point within the nonlocal horizon of $\textbf{X}$. For the rigid body motion given in Eq. \eqref{eq: rigid}, the first integer-order derivative $D^1_{R_j} \Psi_{i}(\textbf{R})$ is
\begin{equation}
D^1_{R_j} \Psi_{i}(\textbf{R})=Q_{ij}
\end{equation}
As discussed in \S\ref{ssec: FCM}, the length-scales $L_{A_j}$ and $L_{B_j}$ in Eq. \eqref{eq: def_grad} are given by $L_{A_{j}}=X_j-X_{A_{j}}$ and $L_{B_{j}}=X_{B_{j}}-X_j$. Using these relations and recurrence definition for the Gamma function $\Gamma(2-\alpha)$, it immediately follows that the deformation gradient tensor is $\tilde{\textbf{F}}_{X_{ij}}^{\Psi}=Q_{ij}$ identically. As discussed previously, this results in null strain measures for rigid body translation and rotations using the fractional-order kinematic relations. Repeating the above procedure, similar arguments can be made for the frame invariance of $\tilde{\textbf{F}}^{\Psi}_x$ and subsequently of $\overset{\alpha}{{\textbf{F}}}$. This establishes the frame-invariance of the model. Note that the above proof holds true independent of the values of $L_{A_j}$ and $L_{B_j}$, and for a more general case $L_{A_j} \neq L_{B_j}$, which occurs in the presence of asymmetric horizons noted at points close to material boundaries and interfaces.

\section*{Appendix 2}
\textbf{Derivation of governing equations: }In the following, we provide the detailed steps for obtaining the strong form of the governing differential equations given in Eq. \eqref{eq: all_governing_equations}. The governing equations are obtained by minimizing the potential energy given in Eq. \eqref{eq: pot_energy}. The minimization of the potential energy results in two independent integrals which are given in Eq.~(\ref{eq: weak_form_eq}) because of the independent nature of the variations $\delta u_0 (x_1)$ and $\delta w_0 (x_1)$. It follows from Eq.~(\ref{eq: weak_form_axial_eq}) that:
\begin{subequations}
\begin{equation}
\label{eq: app_by_parts_1}
    \int_{0}^{L} N(x_1)D_{x_1}^{\alpha}\left[\delta u_0(x_1)\right]\mathrm{d}x_1=-\int_0^L {}^{R-RL}_{x_1-l_B}D_{x_1+l_A}^{\alpha} \left[N(x_1)\right] \delta u_0(x_1) \mathrm{d}x_1+\left[{}_{x_1-l_B}I^{1-\alpha}_{x_1+l_A}\left[N(x_1)\right]\delta u_0\right]\big\vert_{0}^{L}
\end{equation}
In the above equation, the Riesz fractional integral of $N(x_1)$ is given as:
\begin{equation}
\label{eq: reisz integral_def}
    {}_{x_1-l_B}I^{1-\alpha}_{x_1+l_A} \left[N(x_1)\right] =\frac{\Gamma(2-\alpha)}{2}\left( l_B^{\alpha-1} {}_{x_1-l_B}I_{x_1}^{1-\alpha} \left[N(x_1)\right]-l_A^{\alpha-1} {}_{x_1}I_{x_1+l_A}^{1-\alpha} \left[N(x_1)\right]\right)
\end{equation}
\end{subequations}
where ${}_{x_1-l_B}I_{x_1}^{1-\alpha}$ and ${}_{x_1}I_{x_1+l_A}^{1-\alpha}$ are the left and right Riesz integrals \cite{podlubny1998fractional}. The simplification of Eq.~(\ref{eq: weak_form_eq}) to Eq.~(\ref{eq: app_by_parts_1}) is performed using integration by parts and the definitions for the fractional derivatives. The specific steps leading to the above simplification can be found in \cite{patnaik2019FEM}. 
The above discussed simplification extends directly to other integrals in Eq.~(\ref{eq: weak_form_eq}). Combining all the simplifications yields the first variation of the potential energy as:
\begin{equation}
\label{eq: app_total}
\begin{split}
    \delta \Pi&=-\int_0^L \underbrace{\left[{}^{R-RL}_{x_1-l_B}D^{\alpha}_{x_1+l_A}[N(x_1)]+F_a(x_1)\right]}_{\text{Axial Governing Equation}} \delta u_0(x_1) \mathrm{d}x_1 + \underbrace{\left[{}_{x_1-l_B}I^{1-\alpha}_{x_1+l_A}[N(x_1)]\delta u_0(x_1)\right]\Big\vert_{0}^{L}}_{\text{Axial Boundary Conditions}}\\
    & - \int_0^L \underbrace{\left[D^1_{x_1} \left({}^{R-RL}_{x_1-l_B} D^{\alpha}_{x_1+l_A}[M(x_1)]\right)+ {}^{R-RL}_{x_1-l_B}D^{\alpha}_{x_1+l_A} \left[N(x_1)D_{x_1}^{\alpha}[w_0(x_1)]\right] + F_t(x_1)\right]}_{\text{Transverse Governing Equation}} \delta w_0(x_1) \mathrm{d}x_1 \\
    & - \underbrace{\left[{}_{x_1-l_B}I^{1-\alpha}_{x_1+l_A}[M(x_1)]  \delta D^1_{x_1}w_0(x_1) - \left\{{}_{x_1-l_B}^{R-RL}D^{\alpha}_{x_1+l_A}[M(x_1)]+{}_{x_1-l_B}I^{1-\alpha}_{x_1+l_A}[N(x_1)] \right\} \delta w_0(x_1)\right]\Big\vert_{0}^{L}}_{\text{Transverse Boundary Conditions}}
\end{split}
\end{equation}
where $D^1_{x_1}(\cdot)$ denotes the first integer-order derivative with respect to $x_1$. The governing differential equations and the associated boundary conditions given in Eq. \eqref{eq: all_governing_equations} are now obtained from the above Eq.~(\ref{eq: app_total}) by using first principle of variational calculus. Note that, for longitudinal boundaries of the beam, i.e., $x_1=\{0,L\}$, $l_A$ and $l_B$ are equal to zero in the Caputo and RL derivatives defined in Eqs.~(\ref{eq: frac_der_def},\ref{eq: r_rl_frac_der_def}), respectively, and the Riesz integrals given in Eq. \eqref{eq: reisz integral_def}. These limit cases reduces the boundary conditions to the integer-order counterparts as given in Eq.~\eqref{eq: all_governing_equations} \cite{patnaik2019FEM}. 
\section{Acknowledgements}
The following work was supported by the Defense Advanced Research Project Agency (DARPA) under the grant \#D19AP00052, and the National Science Foundation (NSF) under the grant DCSD \#1825837. The content and information presented in this manuscript do not necessarily reflect the position or the policy of the government. The material is approved for public release; distribution is unlimited.

\bibliographystyle{unsrt}
\bibliography{nonlinear_ffem}

\begin{thebibliography}{10}

\bibitem{szabo1994time}
T.~L. Szabo.
\newblock Time domain wave equations for lossy media obeying a frequency power
  law.
\newblock {\em The Journal of the Acoustical Society of America},
  96(1):491--500, 1994.

\bibitem{adbone}
K.~A. Wear.
\newblock A stratified model to predict dispersion in trabecular bone.
\newblock {\em IEEE transactions on ultrasonics, ferroelectrics, and frequency
  control}, 48(4):1079--1083, 2001.

\bibitem{fellah2004verification}
Z.~E.~A. Fellah, S.~Berger, W.~Lauriks, and C.~Depollier.
\newblock Verification of kramers--kronig relationship in porous materials
  having a rigid frame.
\newblock {\em Journal of sound and vibration}, 270(4-5):865--885, 2004.

\bibitem{stulov2016frequency}
A.~Stulov and V.~I. Erofeev.
\newblock Frequency-dependent attenuation and phase velocity dispersion of an
  acoustic wave propagating in the media with damages.
\newblock In {\em Generalized Continua as Models for Classical and Advanced
  Materials}, pages 413--423. Springer, 2016.

\bibitem{kroner1967elasticity}
E~Kr{\"o}ner.
\newblock Elasticity theory of materials with long range cohesive forces.
\newblock {\em International Journal of Solids and Structures}, 3(5):731--742,
  1967.

\bibitem{eringen1972nonlocal}
A.~C. Eringen and D.~G.~B. Edelen.
\newblock On nonlocal elasticity.
\newblock {\em International Journal of Engineering Science}, 10(3):233--248,
  1972.

\bibitem{sumelka2015fractional}
W.~Sumelka, T.~Blaszczyk, and C.~Liebold.
\newblock Fractional euler--bernoulli beams: Theory, numerical study and
  experimental validation.
\newblock {\em European Journal of Mechanics-A/Solids}, 54:243--251, 2015.

\bibitem{civalek2016simple}
{\"O}.~Civalek and C.~Demir.
\newblock A simple mathematical model of microtubules surrounded by an elastic
  matrix by nonlocal finite element method.
\newblock {\em Applied Mathematics and Computation}, 289:335--352, 2016.

\bibitem{rahimi2017linear}
Z.~Rahimi, W.~Sumelka, and X-J. Yang.
\newblock Linear and non-linear free vibration of nano beams based on a new
  fractional non-local theory.
\newblock {\em Engineering Computations}, 34(5):1754--1770, 2017.

\bibitem{peerlings2001critical}
R.~H.~J. Peerlings, M.~G.~D. Geers, R.~De~Borst, and W.~A.~M. Brekelmans.
\newblock A critical comparison of nonlocal and gradient-enhanced softening
  continua.
\newblock {\em International Journal of solids and Structures},
  38(44-45):7723--7746, 2001.

\bibitem{aifantis2003update}
E.~C. Aifantis.
\newblock Update on a class of gradient theories.
\newblock {\em Mechanics of materials}, 35(3-6):259--280, 2003.

\bibitem{guha2015review}
Suman Guha, Sandeep Sangal, and Sumit Basu.
\newblock A review of higher order strain gradient theories of plasticity:
  Origins, thermodynamics and connections with dislocation mechanics.
\newblock {\em Sadhana}, 40(4):1205--1240, 2015.

\bibitem{sidhardh2018element}
S.~Sidhardh and M.~C. Ray.
\newblock Element-free galerkin model of nano-beams considering strain gradient
  elasticity.
\newblock {\em Acta Mechanica}, 229(7):2765--2786, 2018.

\bibitem{sidhardh2019size}
S.~Sidhardh and M.~C. Ray.
\newblock Size-dependent elastic response in functionally graded microbeams
  considering generalized first strain gradient elasticity.
\newblock {\em The Quarterly Journal of Mechanics and Applied Mathematics},
  2019.

\bibitem{polizzotto2001nonlocal}
C.~Polizzotto.
\newblock Nonlocal elasticity and related variational principles.
\newblock {\em International Journal of Solids and Structures},
  38(42-43):7359--7380, 2001.

\bibitem{bavzant2002nonlocal}
Z.~P. Ba{\v{z}}ant and M.~Jir{\'a}sek.
\newblock Nonlocal integral formulations of plasticity and damage: survey of
  progress.
\newblock {\em Journal of Engineering Mechanics}, 128(11):1119--1149, 2002.

\bibitem{sidhardh2018effect}
S.~Sidhardh and M.~C. Ray.
\newblock Effect of nonlocal elasticity on the performance of a flexoelectric
  layer as a distributed actuator of nanobeams.
\newblock {\em International Journal of Mechanics and Materials in Design},
  14(2):297--311, 2018.

\bibitem{yang2010nonlinear}
J.~Yang, L.~L. Ke, and S.~Kitipornchai.
\newblock Nonlinear free vibration of single-walled carbon nanotubes using
  nonlocal timoshenko beam theory.
\newblock {\em Physica E: Low-dimensional Systems and Nanostructures},
  42(5):1727--1735, 2010.

\bibitem{srividhya2018nonlocal}
S.~Srividhya, P.~Raghu, A.~Rajagopal, and J.~N. Reddy.
\newblock Nonlocal nonlinear analysis of functionally graded plates using
  third-order shear deformation theory.
\newblock {\em International Journal of Engineering Science}, 125:1--22, 2018.

\bibitem{reddy2010nonlocal}
J.~N. Reddy.
\newblock Nonlocal nonlinear formulations for bending of classical and shear
  deformation theories of beams and plates.
\newblock {\em International Journal of Engineering Science},
  48(11):1507--1518, 2010.

\bibitem{emam2013general}
S.~A. Emam.
\newblock A general nonlocal nonlinear model for buckling of nanobeams.
\newblock {\em Applied Mathematical Modelling}, 37(10-11):6929--6939, 2013.

\bibitem{lembo2016nonlinear}
M.~Lembo.
\newblock On nonlinear deformations of nonlocal elastic rods.
\newblock {\em International Journal of Solids and Structures}, 90:215--227,
  2016.

\bibitem{srinivasa2013model}
A.~R. Srinivasa and J.~N. Reddy.
\newblock A model for a constrained, finitely deforming, elastic solid with
  rotation gradient dependent strain energy, and its specialization to von
  k{\'a}rm{\'a}n plates and beams.
\newblock {\em Journal of the Mechanics and Physics of Solids}, 61(3):873--885,
  2013.

\bibitem{reddy2014eringen}
J.~N. Reddy and S.~El-Borgi.
\newblock Eringen’s nonlocal theories of beams accounting for moderate
  rotations.
\newblock {\em International Journal of Engineering Science}, 82:159--177,
  2014.

\bibitem{najar2015nonlinear}
F.~Najar, S.~El-Borgi, J.~N. Reddy, and K.~Mrabet.
\newblock Nonlinear nonlocal analysis of electrostatic nanoactuators.
\newblock {\em Composite Structures}, 120:117--128, 2015.

\bibitem{reddy2007nonlocal}
J.~N. Reddy.
\newblock Nonlocal theories for bending, buckling and vibration of beams.
\newblock {\em International Journal of Engineering Science}, 45(2-8):288--307,
  2007.

\bibitem{challamel2008small}
N.~Challamel and C.~M. Wang.
\newblock The small length scale effect for a non-local cantilever beam: a
  paradox solved.
\newblock {\em Nanotechnology}, 19(34):345703, 2008.

\bibitem{khodabakhshi2015unified}
P.~Khodabakhshi and J.~N. Reddy.
\newblock A unified integro-differential nonlocal model.
\newblock {\em International Journal of Engineering Science}, 95:60--75, 2015.

\bibitem{romano2017constitutive}
Giovanni Romano, Raffaele Barretta, Marina Diaco, and Francesco~Marotti
  de~Sciarra.
\newblock Constitutive boundary conditions and paradoxes in nonlocal elastic
  nanobeams.
\newblock {\em International Journal of Mechanical Sciences}, 121:151--156,
  2017.

\bibitem{ansari2015size}
R.~Ansari, R.~Gholami, and H.~Rouhi.
\newblock Size-dependent nonlinear forced vibration analysis of
  magneto-electro-thermo-elastic timoshenko nanobeams based upon the nonlocal
  elasticity theory.
\newblock {\em Composite Structures}, 126:216--226, 2015.

\bibitem{bagley1983theoretical}
R.~L. Bagley and P.~J. Torvik.
\newblock A theoretical basis for the application of fractional calculus to
  viscoelasticity.
\newblock {\em Journal of Rheology}, 27(3):201--210, 1983.

\bibitem{chatterjee2005statistical}
A.~Chatterjee.
\newblock Statistical origins of fractional derivatives in viscoelasticity.
\newblock {\em Journal of Sound and Vibration}, 284(3-5):1239--1245, 2005.

\bibitem{patnaik2020application}
Sansit Patnaik and Fabio Semperlotti.
\newblock Application of variable-and distributed-order fractional operators to
  the dynamic analysis of nonlinear oscillators.
\newblock {\em Nonlinear Dynamics}, pages 1--20, 2020.

\bibitem{mainardi1996fractional}
F.~Mainardi.
\newblock Fractional relaxation-oscillation and fractional diffusion-wave
  phenomena.
\newblock {\em Chaos, Solitons \& Fractals}, 7(9):1461--1477, 1996.

\bibitem{buonocore2018occurrence}
S.~Buonocore, M.~Sen, and F.~Semperlotti.
\newblock Occurrence of anomalous diffusion and non-local response in
  highly-scattering acoustic periodic media.
\newblock {\em New Journal of Physics}, 2019.

\bibitem{hollkamp2019analysis}
John~P Hollkamp, Mihir Sen, and Fabio Semperlotti.
\newblock Analysis of dispersion and propagation properties in a periodic rod
  using a space-fractional wave equation.
\newblock {\em Journal of Sound and Vibration}, 441:204--220, 2019.

\bibitem{hollkamp2018model}
J.~P. Hollkamp, M.~Sen, and F.~Semperlotti.
\newblock Model-order reduction of lumped parameter systems via fractional
  calculus.
\newblock {\em Journal of Sound and Vibration}, 419:526--543, 2018.

\bibitem{hollkamp2020application}
J.~P. Hollkamp and F.~Semperlotti.
\newblock Application of fractional order operators to the simulation of ducts
  with acoustic black hole terminations.
\newblock {\em Journal of Sound and Vibration}, 465:115035, 2020.

\bibitem{carpinteri2011fractional}
A.~Carpinteri, P.~Cornetti, and A.~Sapora.
\newblock A fractional calculus approach to nonlocal elasticity.
\newblock {\em The European Physical Journal Special Topics}, 193(1):193, 2011.

\bibitem{sumelka2014thermoelasticity}
W.~Sumelka.
\newblock Thermoelasticity in the framework of the fractional continuum
  mechanics.
\newblock {\em Journal of Thermal Stresses}, 37(6):678--706, 2014.

\bibitem{sumelka2016fractional}
W.~Sumelka.
\newblock Fractional calculus for continuum mechanics--anisotropic
  non-locality.
\newblock {\em Bulletin of the Polish Academy of Sciences Technical Sciences},
  64(2):361--372, 2016.

\bibitem{laskin2006nonlinear}
Nick Laskin and G~Zaslavsky.
\newblock Nonlinear fractional dynamics on a lattice with long range
  interactions.
\newblock {\em Physica A: Statistical Mechanics and its Applications},
  368(1):38--54, 2006.

\bibitem{di2008long}
M.~Di~Paola and M.~Zingales.
\newblock Long-range cohesive interactions of non-local continuum faced by
  fractional calculus.
\newblock {\em International Journal of Solids and Structures},
  45(21):5642--5659, 2008.

\bibitem{cottone2009fractional}
G.~Cottone, M.~Di~Paola, and M.~Zingales.
\newblock Fractional mechanical model for the dynamics of non-local continuum.
\newblock In {\em Advances in numerical methods}, pages 389--423. Springer,
  2009.

\bibitem{alotta2014finite}
Gioacchino Alotta, Giuseppe Failla, and Massimiliano Zingales.
\newblock Finite element method for a nonlocal timoshenko beam model.
\newblock {\em Finite Elements in Analysis and Design}, 89:77--92, 2014.

\bibitem{alotta2017finite}
G.~Alotta, G.~Failla, and M.~Zingales.
\newblock Finite-element formulation of a nonlocal hereditary fractional-order
  timoshenko beam.
\newblock {\em Journal of Engineering Mechanics}, 143(5):D4015001, 2017.

\bibitem{patnaik2019FEM}
Sansit Patnaik, Sai Sidhardh, and Fabio Semperlotti.
\newblock A {R}itz-based finite element method for a fractional-order boundary
  value problem of nonlocal elasticity.
\newblock {\em International Journal of Solids and Structures}, 2020.

\bibitem{szajek2020selected}
Krzysztof Szajek, Wojciech Sumelka, Tomasz Blaszczyk, and Krzysztof Bekus.
\newblock On selected aspects of space-fractional continuum mechanics model
  approximation.
\newblock {\em International Journal of Mechanical Sciences}, 167:105287, 2020.

\bibitem{patnaik2020fractional}
Sansit Patnaik, Sai Sidhardh, and Fabio Semperlotti.
\newblock Fractional-order models for the static and dynamic analysis of
  nonlocal plates.
\newblock {\em arXiv preprint arXiv:2002.10244}, 2020.

\bibitem{sidhardh2020thermoelastic}
Sai Sidhardh, Sansit Patnaik, and Fabio Semperlotti.
\newblock Thermoelastic response of fractional-order nonlocal and geometrically
  nonlinear beams.
\newblock {\em arXiv preprint arXiv:2003.10215}, 2020.

\bibitem{zheng2010note}
Y.~Zheng, C.~Li, and Z.~Zhao.
\newblock A note on the finite element method for the space-fractional
  advection diffusion equation.
\newblock {\em Computers \& Mathematics with Applications}, 59(5):1718--1726,
  2010.

\bibitem{jin2016petrov}
B.~Jin, R.~Lazarov, and Z.~Zhou.
\newblock A petrov--galerkin finite element method for fractional
  convection-diffusion equations.
\newblock {\em SIAM Journal on Numerical Analysis}, 54(1):481--503, 2016.

\bibitem{reddy2014introduction}
J.~N. Reddy.
\newblock {\em An Introduction to Nonlinear Finite Element Analysis: with
  applications to heat transfer, fluid mechanics, and solid mechanics}.
\newblock OUP Oxford, 2014.

\bibitem{pisano2009nonlocal}
A.~A. Pisano, A.~Sofi, and P.~Fuschi.
\newblock Nonlocal integral elasticity: 2d finite element based solutions.
\newblock {\em International Journal of Solids and Structures},
  46(21):3836--3849, 2009.

\bibitem{norouzzadeh2017finite}
A.~Norouzzadeh and R.~Ansari.
\newblock Finite element analysis of nano-scale timoshenko beams using the
  integral model of nonlocal elasticity.
\newblock {\em Physica E: Low-dimensional Systems and Nanostructures},
  88:194--200, 2017.

\bibitem{zhuang2009numerical}
Pinghui Zhuang, Fawang Liu, Vo~Anh, and Ian Turner.
\newblock Numerical methods for the variable-order fractional
  advection-diffusion equation with a nonlinear source term.
\newblock {\em SIAM Journal on Numerical Analysis}, 47(3):1760--1781, 2009.

\bibitem{zhang2013novel}
H~Zhang, Fawang Liu, Mantha~S Phanikumar, and Mark~M Meerschaert.
\newblock A novel numerical method for the time variable fractional order
  mobile--immobile advection--dispersion model.
\newblock {\em Computers \& Mathematics with Applications}, 66(5):693--701,
  2013.

\bibitem{fernandez2016bending}
J.~Fern{\'a}ndez-S{\'a}ez, R.~Zaera, J.~A. Loya, and J.~N. Reddy.
\newblock Bending of euler--bernoulli beams using eringen’s integral
  formulation: a paradox resolved.
\newblock {\em International Journal of Engineering Science}, 99:107--116,
  2016.

\bibitem{sumelka2015non}
W~Sumelka.
\newblock Non-local kirchhoff--love plates in terms of fractional calculus.
\newblock {\em Archives of Civil and Mechanical Engineering}, 15(1):231--242,
  2015.

\bibitem{podlubny1998fractional}
I.~Podlubny.
\newblock {\em Fractional differential equations: an introduction to fractional
  derivatives, fractional differential equations, to methods of their solution
  and some of their applications}, volume 198.
\newblock Elsevier, 1998.

\end{thebibliography}

\end{document}